\documentclass[12pt]{article}

\usepackage{graphicx}
\usepackage{subfigure}
\usepackage{psfrag}
\usepackage{rotating}

\usepackage{amsfonts}
\usepackage{amssymb}
\usepackage{amsmath}
\usepackage{amsthm}

\newcommand{\beq}{\begin{equation}}
\newcommand{\eeq}{\end{equation}}
\newcommand{\gf}{generating function}

\theoremstyle{plain}
\newtheorem{theo}{\bf Theorem}[section]
\newtheorem{propo}[theo]{\bf Proposition}

\newtheorem{coro}[theo]{\bf Corollary}
\newtheorem{lem}[theo]{\bf Lemma}
\newtheorem{defin}[theo]{\bf Definition}

\theoremstyle{definition}
\newtheorem{exam}[theo]{\bf Example}

\newcommand{\Ref}[1]{(\ref{#1})}

\newcommand{\field}[1]{\mathbb{#1}}

\newcommand{\fR}{\field{R}}

\def\qed{$\hfill{\vrule height 3pt width 5pt depth 2pt}$}

\def\ssA{{\mathbf{\mathsf{A}}}}

\def\ssT{{\mathbf{\mathsf{T}}}}

\def\Mt{{\overline M}}

\def\Mb{{\underline M}}

\def\Nt{{\overline N}}

\def\yt{{\overline y}}
\def\wt{{\overline w}}
\def\zt{{\overline z}}
\def\Nb{{\underline N}}

\def\yb{{\underline y}}
\def\wb{{\underline w}}
\def\zb{{\underline z}}
\def\St{{\overline S}}

\def\Sb{{\underline S}}

\def\Ht{{\overline H}}

\def\Hb{{\underline H}}

\def\bM{\mathbf{M}}
\def\bm{{\boldsymbol m}}
\def\bn{{\boldsymbol n}}
\def\bMt{\mathbf{{\overline M}}}

\def\bMb{\mathbf{{\underline M}}}

\def\bNt{\mathbf{{\overline N}}}

\def\byt{\boldsymbol{{\overline y}}}
\def\byb{\boldsymbol{{\underline y}}}
\def\bwt{\boldsymbol{{\overline w}}}
\def\bwb{\boldsymbol{{\underline w}}}
\def\bzb{\boldsymbol{{\underline z}}}
\def\bzt{\boldsymbol{{\overline z}}}
\def\bNb{\mathbf{{\underline N}}}

\def\bSt{\mathbf{{\overline S}}}

\def\bSb{\mathbf{{\underline S}}}

\def\bHt{\mathbf{{\overline H}}}

\def\bHb{\mathbf{{\underline H}}}

\def\bzero{{\mathbf 0}}
\def\bone{{\mathbf 1}}

\def\bu{{\boldsymbol u}}

\def\bv{{\boldsymbol v}}

\def\bx{{\boldsymbol x}}

\def\by{{\boldsymbol y}}

\def\bz{{\boldsymbol z}}

\def\balpha{{\boldsymbol{\alpha}}}
\def\bbeta{{\boldsymbol{\beta}}}
\def\bgamma{{\boldsymbol{\gamma}}}
\def\bdelta{{\boldsymbol{\delta}}}
\def\btau{{\boldsymbol{\tau}}}
\def\bepsilon{{\boldsymbol{\epsilon}}}
\def\bPhi{{\boldsymbol{\Phi}}}
\def\bphi{{\boldsymbol{\phi}}}

\def\ap#1#2#3 {{\em Ann. Phys. (NY),} {\bf#1} (19#2) #3}
\def\apj#1#2#3 {{\em Astrophys. J.,} {\bf#1} (19#2) #3}
\def\apjl#1#2#3 {{\em Astrophys. J. Lett.,} {\bf#1} (19#2) #3}
\def\app#1#2#3 {{\em Acta. Phys. Pol.,} {\bf#1} (19#2) #3}
\def\ar#1#2#3 {{\em Ann. Rev. Nucl. Part. Sci.,} {\bf#1} (19#2) #3}
\def\cmp#1#2#3 {{\em Comm. Math. Phys.,} {\bf#1} (19#2) #3}
\def\cpc#1#2#3 {{\em Computer Phys. Comm.,} {\bf#1} (19#2) #3}
\def\err#1#2#3 {{\it Erratum} {\bf#1} (19#2) #3}
\def\ib#1#2#3 {{\it ibid.} {\bf#1} (19#2) #3}
\def\jmp#1#2#3 {{\em J. Math. Phys.,} {\bf#1} (19#2) #3}
\def\ijmp#1#2#3 {{\em Int. J. Mod. Phys.,} {\bf#1} (19#2) #3}
\def\jetp#1#2#3 {{\em JETP Lett.,} {\bf#1} (19#2) #3}
\def\jpa#1#2#3 {{\em J. Phys. A: Math. Gen.,} {\bf#1} (19#2) #3}
\def\jpg#1#2#3 {{\em J. Phys. G.,} {\bf#1} (19#2) #3}
\def\jsp#1#2#3 {{\em J. Stat. Phys.,} {\bf#1} (19#2) #3}
\def\mpl#1#2#3 {{\em Mod. Phys. Lett.,} {\bf#1} (19#2) #3}
\def\nat#1#2#3 {{\em Nature (London),} {\bf#1} (19#2) #3}
\def\nc#1#2#3 {{\em Nuovo Cim.,} {\bf#1} (19#2) #3}
\def\nim#1#2#3 {{\em Nucl. Instr. Meth.,} {\bf#1} (19#2) #3}
\def\np#1#2#3 {{\em Nucl. Phys.,} {\bf#1} (19#2) #3}
\def\pcps#1#2#3 {{\em Proc. Cam. Phil. Soc.,} {\bf#1} (#2) #3}
\def\pl#1#2#3 {{\em Phys. Lett.,} {\bf#1} (19#2) #3}
\def\prep#1#2#3 {{\em Phys. Rep.,} {\bf#1} (19#2) #3}
\def\prev#1#2#3 {{\em Phys. Rev.,} {\bf#1} (19#2) #3}
\def\prl#1#2#3 {{\em Phys. Rev. Lett.,} {\bf#1} (19#2) #3}
\def\prs#1#2#3 {{\em Proc. Roy. Soc.,} {\bf#1} (19#2) #3}
\def\ptp#1#2#3 {{\em Prog. Th. Phys.,} {\bf#1} (19#2) #3}
\def\ps#1#2#3 {{\em Physica Scripta,} {\bf#1} (19#2) #3}
\def\rmp#1#2#3 {{\em Rev. Mod. Phys.,} {\bf#1} (19#2) #3}
\def\rpp#1#2#3 {{\em Rep. Prog. Phys.,} {\bf#1} (19#2) #3}
\def\sjnp#1#2#3 {{\em Sov. J. Nucl. Phys.,} {\bf#1} (19#2) #3}
\def\spj#1#2#3 {{\em Sov. Phys. JETP,} {\bf#1} (19#2) #3}
\def\spu#1#2#3 {{\em Sov. Phys. Usp.,} {\bf#1} (19#2) #3}
\def\zp#1#2#3 {{\em Zeit. Phys.,} {\bf#1} (19#2) #3}

\begin{document}

\begin{titlepage}
\begin{center}
{\Large \bf Laboratoire Bordelais de Recherche en Informatique \\
Universit{\'e} Bordeaux I \\[0.1in]}
and \\[0.1in]
{\Large \bf Department of Mathematics and Statistics \\
The University of Melbourne \\[0.1in]}
\end{center}
\begin{flushright}
June 1999 \\[0.2in]
\end{flushright}
\begin{center}
{\Large \bf  Inversion Relations, Reciprocity and Polyominoes \\[0.2in]}
{M. Bousquet-M{\'e}lou$^{(1)}$,
A.J. Guttmann$^{(2)}$, W.P. Orrick$^{(2)}$ and A. Rechnitzer$^{(2)}$\\[0.1in]}
$^{(1)}$LaBRI, CNRS \\
Universit{\'e} Bordeaux I \\
351 cours de la Lib{\'e}ration \\
33405 Talence Cedex \\
France \\[0.1in]
$^{(2)}$Department of Mathematics and Statistics\\
The University of Melbourne\\
Parkville, Vic. 3052 \\
Australia \\[0.1in]
\end{center}

\begin{abstract}
We derive self-reciprocity properties for a number of polyomino \gf s,
including several families of 
column-convex polygons, three-choice polygons and staircase polygons
with a staircase hole.  In so doing,
we establish a connection between the reciprocity results known
to combinatorialists and the inversion relations used by physicists to
solve models in statistical mechanics.  For several classes of convex
polygons, the inversion (reciprocity) relation, augmented by
certain symmetry and analyticity properties, completely determines the
anisotropic perimeter generating function.


\end{abstract}

\end{titlepage}

\section{Introduction}
Symmetries are among the most important guiding principles in all of
physics and mathematics.  It often happens that a problem may be solved
by symmetry considerations alone, and even if not,
understanding the symmetries of the solution can greatly reduce the amount
of work needed to find it.  
We study here a symmetry of functions which is known as  ``self-reciprocity''
to combinatorialists and which is referred to as ``inversion relations''
in lattice statistical mechanics.

Our focus will be on polyomino enumeration problems which are of
interest
 in both combinatorics and physics.
We shall demonstrate that one can find examples of functional symmetry
in the resulting \gf s. 

The inversion relation rose to prominence in statistical mechanics in
the early 1980s as the most direct path to the
solution of many integrable
models~\cite{strog,baxb,baxbbis} and was soon realized to be commonplace in
both solved and unsolved models~\cite{baxb,baxbbis,maila}.  Let $G(\bx)$ be a
thermodynamic
quantity which depends on a collection of parameters, $\bx$.
An inversion relation is a functional equation
\begin{equation}
G(\bx) \pm \bx^{\balpha}G(\bphi(\bx)) = \psi(\bx)
\label{invrel}
\end{equation}
where $\bphi$ and $\psi$ are known functions of $\bx$.  Typically $\bphi$
involves taking reciprocals of one or more components of $\bx$.  The
inversion relation tightly constrains the function $G$.  For some
two-dimensional models a pair of additional conditions holds: that
$G$ is symmetric under exchange of horizontal and vertical, and that
$G$ is an analytic function of its arguments.
Very often, the three constraints taken together uniquely
determine the function $G$.

  In 1974  Stanley presented a
general framework for reciprocity results.  He
established several powerful general conditions under which a
generating function will be self-reciprocal~\cite{stanb}.  The
language and notation of Stanley~\cite{stanb,stane} will be used
throughout this paper.

\begin{defin} Let $H(y_1, \ldots , y_n)$ be a rational function in the
variables $y_i$, with coefficients in $\fR$. We say
that $H$ is {\em self-reciprocal} if there exists an $n$-tuple of
integers $(\beta_1 , \ldots , \beta _n)$ such that
\begin{equation}
H(1/y_1, \ldots , 1/y_n) = \pm y_1^{\beta_1} \ldots  y_n^{\beta_n}
H(y_1, \ldots , y_n).
\label{selfrec}
\end{equation}
\label{selfrecdef}
\end{defin}
\noindent
In what follows, we write
$\by^\bbeta\equiv y_1^{\beta_1}\ldots y_n^{\beta_n}$ and
$1/{\mathbf y}\equiv (1/y_1,\ldots,1/y_n)$.  Thus eqn.~(\ref{selfrec}) may
be concisely expressed as
$$
H(1/\by)=\pm \by^\bbeta H(\by).
$$
Note that a rational function is self-reciprocal if and only if both
its numerator and denominator are so, and that the self-reciprocity of a
polynomial amounts to a certain symmetry in its coefficients. Some
explicit examples are given in Subsection \ref{subsect:examples}.

Let us now demonstrate the
relationship between self-reciprocity and inversion relations.
Consider the multivariable generating function
\begin{align}
G(\bx,\by)&=\sum_{\bm,\bn} C(\bm,\bn) x_1^{m_1}\ldots x_j^{m_j} y_1^{n_1}\ldots
 y_k^{n_k} \notag\\
&\equiv \sum_{\bm,\bn} C(\bm,\bn) \bx^{\bm} \by^{\bn}
\label{genfun}
\end{align}
where $\bm=(m_1,\ldots,m_j)$, $\bx=(x_1,\ldots,x_j)$ and similarly for $\bn$
and $\by$.  The summation is over $(j+k)$-tuples of nonnegative integers
representing the objects being enumerated.
Performing the summation over $\bn$, we reexpress
eqn.~(\ref{genfun}) in terms of partial generating functions, $H_{\bm}(\by)$,
\begin{equation}
G(\bx,\by) = \sum_\bm H_\bm(\by) \bx^{\bm}.
\label{pargenfun}
\end{equation}
Now suppose that the partial generating functions are self-reciprocal,
\begin{equation}
H_{\bm}(1/\by) = \pm \bepsilon^\bm \by^{\bbeta(\bm)} H_\bm(\by),
\label{Hrec}
\end{equation}
where $\bepsilon$ is a $j$-tuple of elements in the set $\{-1,1\}$ which
characterizes the dependence of the sign on $\bm$, and where
$\bbeta(\bm)$ depends linearly on $\bm$:
\begin{equation}
\bbeta(\bm)=\ssA\bm+\balpha.
\label{betadef}
\end{equation}
Here, $\ssA=(a_{\ell,i})_{\ell,i}$ is a $k\times j$ matrix of integers and $\balpha$ is a $k$-tuple of
integers.
We can then write
\begin{equation}
G(\bx,\by) \mp \by^{-\balpha} G(\bepsilon\bx \by^{-\ssA},1/\by)=0
\label{Ginversion}
\end{equation}
where $\bepsilon\bx\by^{-\ssA}$ is the $j$-tuple whose $i^{\text{th}}$ entry is
$\epsilon_i x_i \prod_{\ell=1}^k y_\ell^{-a_{\ell,i}}$.
This is clearly a special case of the inversion relation~(\ref{invrel}).
A few comments are in order:

\begin{itemize}
\item The right hand side of~(\ref{Ginversion})
is zero, but in the more general situation some of the partial generating
functions, $H_\bm(\by)$ will fail to be self-reciprocal for certain
choices of $\bm$.  If we are fortunate, this will be a small, finite or
otherwise controllable set of cases, and we will be able to compute
the correction term we need to add to the right hand side explicitly.
For many examples in statistical mechanics, this correction term depends
on $\bx$ but not on $\by$.

\item In all of the cases we shall see below, the denominators of our
rational
functions will be a product of terms $(1- \by^{\balpha_j})$, which
are self-reciprocal.  
Stanley has proved that
this denominator form always holds for certain  classes of
problems (see Theorem~4.6.11 of ref.~\cite{stane}).

\item It might be asked which of the concepts, inversion or self-reciprocity,
is the more general.  On one hand, in the derivation of~(\ref{Ginversion})
the dependence of the exponent $\bbeta$ on $\bm$ was assumed to be
linear, which may not always hold, implying that reciprocity is more
fundamental.  On the other hand, the function $\bphi$ 
occurring in~(\ref{invrel}) 
may in principle
be more complicated than $\bx\rightarrow \bepsilon\bx\by^{-\ssA}$,
$\by\rightarrow 1/\by$. In this case, the partial generating functions
might not be self-reciprocal.  
An example is provided in
Section~\ref{sect:statmech} by  the Potts model, but in the polyomino
 examples considered in this
paper, this situation does not arise.
\end{itemize}

We now present a nonexhaustive list of recipes for finding and proving
reciprocity results and inversion relations.
\begin{enumerate}
\item If the generating function (or thermodynamic quantity) is known in
closed form, an inversion relation can be demonstrated directly.  As an
example, we treat the anisotropic perimeter generating function for directed
convex polygons  in this manner in Section~\ref{sect:polyinversion}.
\label{item:exact}
\item For statistical mechanics models which admit a formulation in
terms of a family of commuting transfer matrices, a transformation
of parameters can often be found which inverts the transfer matrix.
The commutativity property then allows the inversion relation to
be derived.  We review this in detail in Section~\ref{sect:statmech},
with the two-dimensional, zero-field Ising model as primary example.
\label{item:tm}
\item In non-integrable models, the transfer matrix will still be
invertible and may suggest a possible inversion relation,
but the required analyticity property is lacking.
Nevertheless, the suggested inversion relation can often be verified
by inspection of the partial generating functions up to some
finite order in the
low-temperature expansion~(\ref{pargenfun}).  The $q>2$ Potts model
inversion relation
discussed in Section~\ref{sect:statmech} was derived this way in
ref.~\cite{jmb}.
Some of the new results
reported in the present paper were initially discovered by this method before
being rederived by one of the other methods.
\label{item:series}
\item The ``Temperley methodology''~\cite{bousc} can be used
to obtain very general reciprocity results for many classes of column-convex
polygons.  The first step is to derive a functional equation for the
generating function  which can be interpreted as the gluing of
an additional column onto the graph.
Step two is to show by induction that appending an
additional column preserves self-reciprocity. 
This is detailed in Section~\ref{sect:temperley}.
\label{item:temperley}
\item If the problem can be posed as a system of linear diophantine
equations, whose solutions are subject to certain types of constraints, we may
apply self-reciprocity theorems due to Stanley~\cite{stanb}.  We have
so far succeeded in applying this method only to families of
directed polyominoes (Section \ref{sect:stanley}),
but it enables us to treat problems which are impossible, or at least
extremely cumbersome, by the method of functional equations.
\label{item:stanley}
\item For combinatorial objects with a rational \gf \ of denominator
$\prod_j (1-\by ^{\balpha_j})$, one can try to explain
 self-reciprocity -- {\em i.e.}, the symmetry of the numerator -- 
by interpreting the numerator combinatorially. 
This has been done by  F{\'e}dou for a family of
objects related to (but distinct from) staircase polygons~\cite{fedo}.
\end{enumerate}

In Section~\ref{sect:statmech} we review the motivation for looking
at inversion relations in statistical mechanics and 
describe the methods used to obtain them.
This will be useful for
making comparisons with the results obtained later, and for suggesting
applications and generalizations of the inversion relations.  
In Section~\ref{sect:polyinversion}, we present examples of
reciprocity results and inversion 
relations for polyominoes, and summarize our main new results.
The technical heart of the paper consists of Section~\ref{sect:temperley}
on the Temperley methodology, and Section~\ref{sect:stanley} on the
application of Stanley's results to polyominoes.

\section{Inversion relations in statistical mechanics}\label{sect:statmech}
The first use of the inversion relation in statistical mechanics was the
solution by Stroganov of certain two dimensional vertex models on the
square lattice~\cite{strog}.  Generalizations of Stroganov's models were
later solved by the same means by Schultz~\cite{schu}.
Shortly after Stroganov, Baxter used a similar method to solve
the hard hexagon model~\cite{baxa} and recognized its broad applicability,
giving the eight-vertex and Ising models as examples~\cite{baxb}.
Subsequently, a number of authors pointed out that many known
solutions to problems in two-dimensional statistical mechanics can
be derived easily using the inversion relation method.  Among these
were Shankar~\cite{shan}, Baxter~\cite{baxc} and Pokrovsky and
Bashilov~\cite{pb}.

It is noteworthy that inversion relations hold also for
models that have not been solved.  Prominent among such models are the
two-dimensional Ising model in a magnetic field whose inversion relation
was found by Baxter~\cite{baxb}, and the three-dimensional Ising
model and noncritical $q$-state Potts model, both of whose inversion
relations were found by Jaekel and Maillard~\cite{jma,jmb}.
What generally distinguishes solved and unsolved models is the
growth rate in the number of poles arising in the partial generating
functions in the expansion~(\ref{pargenfun}), as a function of order.
Roughly speaking, a
more complicated pole structure implies that the number of parameters
needed to specify a given partial generating function is
greater, and makes it less likely that an inversion relation can
completely determine all of them.  Nevertheless, inversion
relations are still invaluable in the study of such problems, not least
because they provide an independent check on series data.

The fundamental problem of statistical mechanics is to calculate the
partition function.  Here we consider vertex models defined
on a square lattice with each bond colored with one of $r$ possible colors.
Each lattice site makes
a contribution to the energy of the system which depends on the colors
of the adjacent bonds.
This defines an $r^4$-vertex model if all possible colorings are permitted.

Stroganov computed the partition function per
site in the thermodynamic limit of several 16- and 81-vertex models.
Consider first a finite lattice (on the 
torus) of
$N$ rows and $M$ columns.  The partition function can be expressed
in terms of the transfer matrix $\ssT_M$ as:
\begin{equation}
Z_{M,N}={\mathrm{Tr}}\left[\left(\ssT_M\right)^N\right]
\label{pf}
\end{equation}
(see \cite{baxd,stane}). Here, $\ssT_M$ is the $r^M\times r^M$ matrix
whose $i,j$th entry is the 
contribution to $Z_{M,N}$ of a single row of $M$ sites connected to the row
below by a set of vertical bonds in configuration $i$ and to the row above by
a set of vertical bonds in configuration $j$.  It depends on the temperature,
$T$, and on $r^4$ parameters specifying the vertex energies. 
In the thermodynamic limit, the partition function per site is given by
\begin{equation}
\kappa = \lim_{M, N\rightarrow\infty} (Z_{M,N})^{1/MN} =\lim_{M\rightarrow\infty} \left(\lambda_M
\right)^{1/M}
\label{pfpersite}
\end{equation}
where $\lambda_M$ is the largest eigenvalue of $\ssT_M$, assumed to be
nondegenerate.

For simplicity let us consider a family of models whose
vertex energies are functions of a single parameter, $b$.
The models solved by Stroganov are
integrable by virtue of the commutativity of the transfer matrices at
different values of this parameter.  This implies that the transfer
matrix eigenvectors are common to all members of the family, and that
the $b$ dependence is only in the eigenvalues.   For this reason $b$
is often called the {\em spectral parameter}.
The key observation is that the
inverse of the transfer matrix in these models
is itself a member
of the commuting family, up to a scale factor
\begin{equation}
\left[ \ssT_M(b)\right] ^{-1} = \psi(b)^{-M} \ {\ssT_M}\left(\phi(b)\right).
\label{Tinv}
\end{equation}
Acting on the eigenvector corresponding to $\lambda_M(b)$
with both sides of eqn.~(\ref{Tinv}) yields the functional equation
\begin{equation}
\kappa(b) \kappa\left(\phi(b)\right) = \psi(b).
\label{func}
\end{equation}
It is the commutativity of the transfer matrices for all values of $b$
that allows the analytical continuation of the function $\kappa$ from
$b$ to $\phi(b)$.
With knowledge of the functions $\psi(b)$ and $\phi(b)$ and using
the analyticity of $\kappa(b)$, Stroganov finds a unique solution, thereby
reproducing Baxter's results for the symmetric eight vertex and homogeneous
ferroelectric models, and obtaining the result for a certain 81-vertex
model~\cite{strog}.

As an illustrative example, we review here the derivation by
Baxter~\cite{baxb} of Onsager's expression for the partition function of the
two-dimensional zero-field Ising model~\cite{onsa}.
Let the square lattice be drawn
at $45^\circ$ to the horizontal and let the couplings between
nearest neighbors along the two lattice directions be
$J$ and $J'.$  Define low temperature variables
\begin{equation}
x = e^{-2K}, \quad y = e^{-2K'}
\ \ \ \ \hbox{with}\ \ \ \ 
K=J/k_BT, \quad K'=J'/k_BT.
\label{Kdef}
\end{equation}
Transfer matrices for different choices of parameters will commute
provided they have the same value of $k=(\sinh 2K \sinh 2K')^{-1}$. 
%
The transformation 
which inverts the transfer matrix is
%
%
\begin{equation} K
 \rightarrow K + \frac{i\pi}{2},   \quad K'  \rightarrow -K',
\label{invcouplings} 
\end{equation}
which does not modify the value of $k$.
Define the reduced partition function per site by
\begin{equation}
\Lambda(x,y)=\exp(-K-K') \kappa(K,K').
\label{redpf}
\end{equation}
Then $\Lambda(x,y)$ obeys the inversion relation
\begin{equation}
\Lambda(x,y) \Lambda(-x,1/y) = 1-x^2.
\label{isinginv}
\end{equation}
Note that $\log\Lambda(x,y)-\frac{1}{2}\log(1-x^2)$ has an inversion relation
of precisely the form~(\ref{Ginversion}).

By the symmetry of the model, we have
\begin{equation}
\Lambda(x,y)=\Lambda(y,x).
\label{isingsym}
\end{equation}
Inspection of the low temperature expansion leads us to conjecture the form
\begin{equation}
\Lambda(x,y) = 1 + \sum_{m\ge 1} \frac{P_m(y^2)}{(1-y^2)^{2m-1}}x^{2m}.
\label{isingform}
\end{equation}
That the coefficient of $x^{2m}$ is a rational function of $y^2$ is apparent
from the nature of the low temperature expansion, but that the denominator
has such a simple form is not expected on general grounds.
Presumably it is a consequence of the condition of commuting transfer
matrices.  Here we take it as a hypothesis.
Then Baxter has shown that the inversion relation~(\ref{isinginv}),
symmetry~(\ref{isingsym}) and the denominator form~(\ref{isingform})
determine $\Lambda(x,y)$ completely.  We present his argument in
Section~\ref{subsect:computegf} where we use it in the context of
polygon enumeration.

Up till now we have been assuming integrability and
in particular we have relied on the property that the transfer matrix
and its inverse are both members of some one-parameter commuting
family.  What about models for which this property doesn't hold?
Since analyticity of $\kappa(b)$ breaks down, the step~(\ref{func}) in
the above derivation is
no longer valid.   However, it is still possible to obtain an inversion
relation by direct analysis of the low-temperature  expansion
of the partition function to some finite order.
As an example, it was shown in ref.~\cite{jmb} that the logarithm of the
reduced partition function per site, $G(x,y)=\ln\Lambda(x,y)$, of the
$q$-state Potts model satisfies the inversion relation
\begin{equation}
G(x,y) + G\left(-\frac{x}{1+(q-2)x},\frac{1}{y}\right) =
 \ln\left(\frac{(1-x)(1+(q-1)x)}{1+(q-2)x}\right).
\label{Pottsinv}
\end{equation}
When $q=2$ this reduces to the Ising model inversion relation~(\ref{isinginv}).
The inversion relations we will be considering in the remainder of the
paper are derived by analysis of the generating function (analogous to
the low temperature expansion) and do not depend on the models being
integrable.

An additional new feature is seen in this Potts model example.
Neglecting for the moment the nonzero right-hand-side of~(\ref{Pottsinv}),
which can
be eliminated by a suitable redefinition of $G(x,y)$, we notice that
when $q>2$
there is no longer an order-by-order cancellation of the partial generating
functions as defined in~(\ref{pargenfun}), but rather cancellation of
combinations of partial generating functions of different orders.
However, we may convert
to self-reciprocal form by defining
\begin{equation}
G'(x,y)=G\left(\frac{x}{1-(q-2)x/2},y\right)
\label{changeofvars}
\end{equation}
under which the inversion relation becomes
\begin{equation}
G'(x,y)+G'(-x,1/y) = \ln\frac{1-q^2x^2/4}{1-(q-2)^2x^2/4}.
\label{newPottsinv}
\end{equation}
In the cases we will look at in this paper, the partial generating functions
turn out to be self-reciprocal in the natural variables of the
problem.
We have not investigated the existence of inversion relations
involving more complicated changes of variables.
%
%
\section{Polyomino enumeration and self-reciprocity} \label{sect:polyinversion}
%
\subsection{Definitions}
The constructions we will consider are defined on the square lattice.  All
are defined only up to translation on the lattice.
Starting
at a lattice site and moving to one of the four nearest neighbors constitutes a
{\it step} which we may identify with the edge connecting the sites.  A
connected sequence of steps is a {\it path} or {\it walk}.
If no lattice site in the path occurs
more than once, the path is {\it self-avoiding}.  If a path returns to
its starting site in the final step, and otherwise does not intersect
itself, the
result is a {\it self-avoiding polygon}.  The number of steps taken is
the {\it perimeter} of the polygon; the number of steps taken in the
vertical direction is the {\it vertical perimeter.} The {\it horizontal
perimeter} is defined similarly.  The {\em area} is the number of cells of the
lattice enclosed by the polygon.

\begin{figure}
  \centering
  \subfigure[Ferrers graph]{
  \begin{minipage}[b]{0.3\textwidth}
    \centering
    \includegraphics{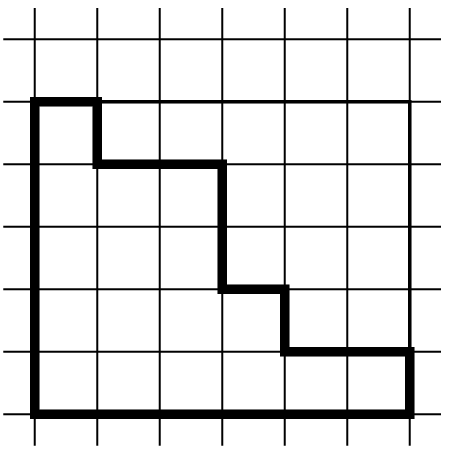}
  \end{minipage}}%
  \subfigure[Stack polygon (horizontal)]{
  \begin{minipage}[b]{0.3\textwidth}
    \centering
    \includegraphics{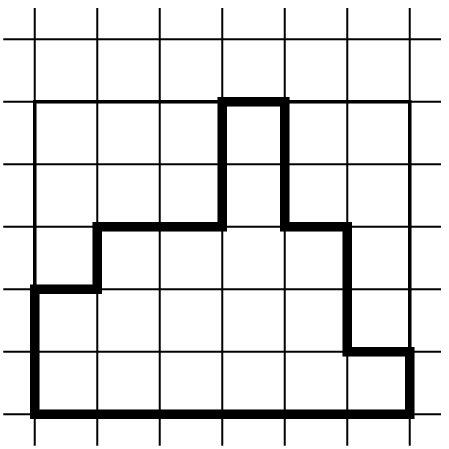}
    \label{subfig:stackv}
  \end{minipage}}%
  \subfigure[Stack polygon (vertical)]{
  \begin{minipage}[b]{0.3\textwidth}
    \centering
    \includegraphics{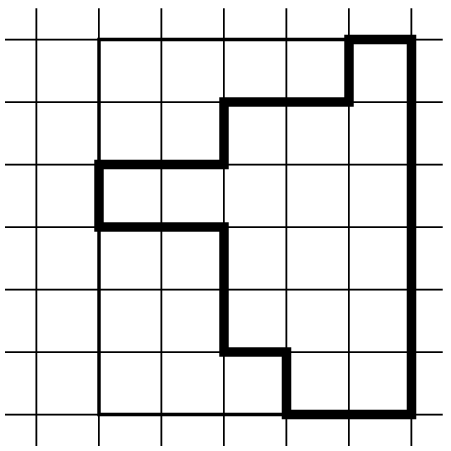}
    \label{subfig:stackh}
  \end{minipage}}\\[20pt]
  \subfigure[Staircase polygon]{
  \begin{minipage}[b]{0.3\textwidth}
    \centering
    \includegraphics{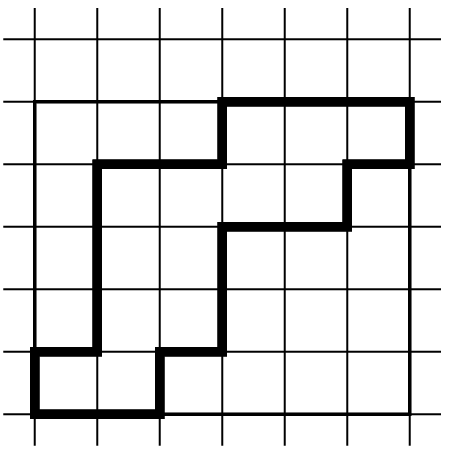}
  \end{minipage}}%
  \subfigure[Directed convex polygon]{
  \begin{minipage}[b]{0.3\textwidth}
    \centering
    \includegraphics{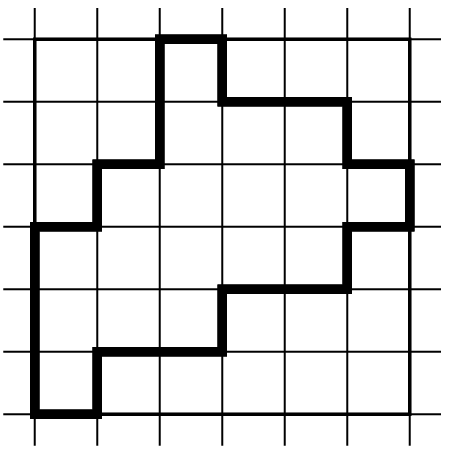}
  \end{minipage}}%
  \subfigure[Convex polygon]{
  \begin{minipage}[b]{0.3\textwidth}
    \centering
    \includegraphics{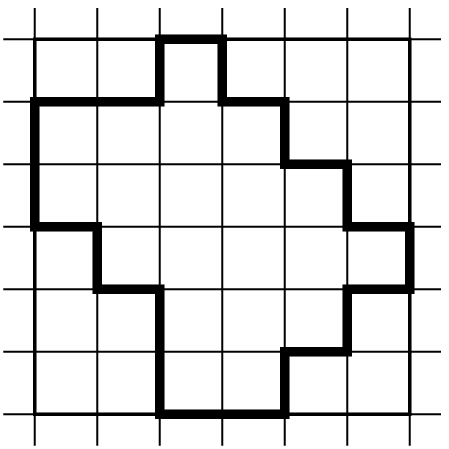}
  \end{minipage}}
  \caption{Classes of convex polygons}
  \label{fig:convexpolygons}
\end{figure}
Enumerating self-avoiding polygons
according to perimeter or area is an unsolved problem.
However, progress has
been made in enumerating certain subclasses
of self-avoiding polygons.  {\em Rectangles} coincide with the rectangles
of ordinary geometry whose vertices are lattice points and whose edges
lie along lattice directions.
A rectangle which contains a given polygon, {\em i.e.}, all steps of the
polygon lie inside or on the rectangle, is a {\it bounding
rectangle} for that polygon.  The smallest such rectangle is the
{\it minimal bounding rectangle.}  A polygon whose perimeter
equals that of its minimal bounding rectangle is {\it convex}.  If a
convex polygon contains at least one of the corners of its minimal
bounding rectangle (for concreteness say the south-west corner) then it
is a {\it directed convex polygon.}  If it contains also the north-east
corner, it is a {\it staircase polygon,} so called because it is bounded
above and below by two staircase-like or {\em directed}
paths.  On the other hand, if it contains two adjacent corners, say the
southwest and southeast (northeast and southeast) then it is a {\em stack
polygon} with horizontal (vertical) orientation.  If it contains three
corners, then it is a {\em Ferrers graph.}
Representative examples of different classes of convex polygons are shown in
Figure~\ref{fig:convexpolygons}.

\begin{figure}
  \centering
  \subfigure[Column-convex polygon]{
  \begin{minipage}[b]{0.3\textwidth}
    \centering
    \includegraphics{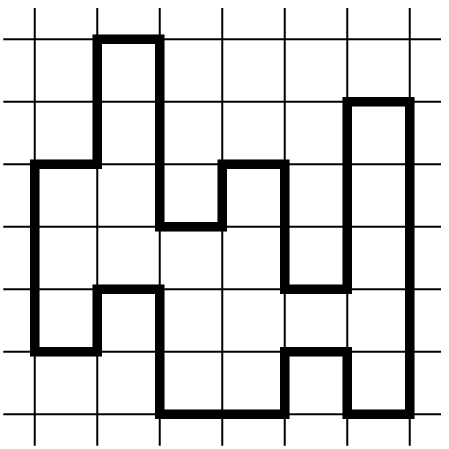}
  \end{minipage}}%
  \subfigure[Bar-graph]{
  \begin{minipage}[b]{0.3\textwidth}
    \centering
    \includegraphics{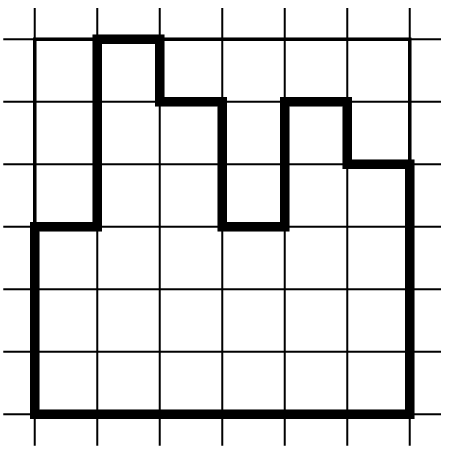}
  \end{minipage}}%
  \subfigure[Directed column-convex polygon]{
  \begin{minipage}[b]{0.3\textwidth}
    \centering
    \includegraphics{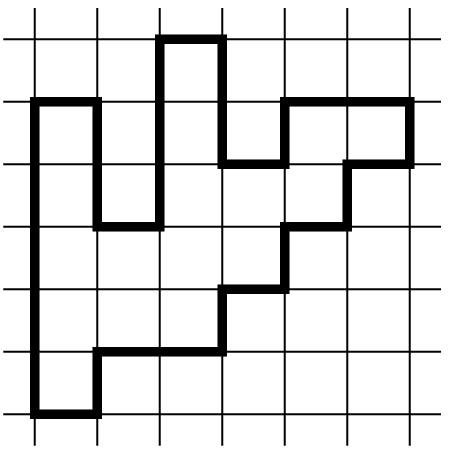}
  \end{minipage}}
  \caption{Classes of column-convex polygons}
  \label{fig:columnconvexpolygons}
\end{figure}
One way to obtain non-convex polygons is to relax
the convexity condition along one direction only.
A self-avoiding polygon is {\it column-convex\/} if 
the intersection of any vertical line with the polygon has at most two
connected components.
{\it Row-convex} polygons are
similarly defined.  The set of convex polygons is the intersection of the
sets of row- and column-convex polygons.  Subclasses of column-convex
polygons include the
{\em bar-graphs} which contain the bottom edge of the
minimal bounding rectangle, and {\em directed column-convex polygons} whose
bottom edge is a directed path.  Some examples are shown in
Figure~\ref{fig:columnconvexpolygons}.

\begin{figure}
  \centering
%
  \subfigure[Three-choice polygon]{
  \begin{minipage}[b]{0.3\textwidth}
    \centering
    \includegraphics{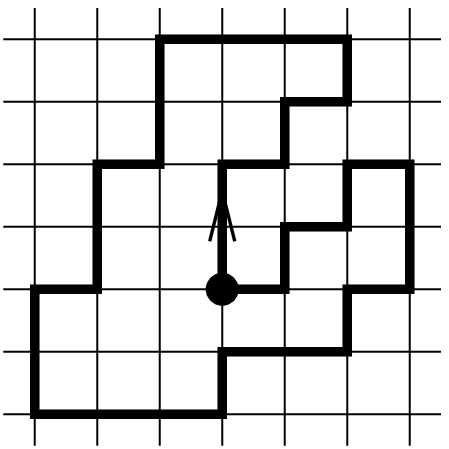}
    \label{subfig:thchoice}
  \end{minipage}}%
  \subfigure[Staircase polygon with a staircase hole]{
  \begin{minipage}[b]{0.3\textwidth}
    \centering
    \includegraphics{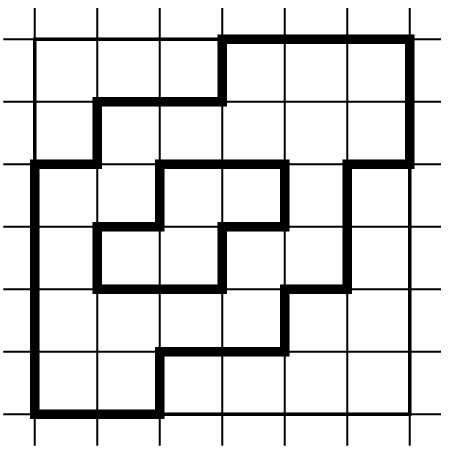}
    \label{subfig:spsh}
  \end{minipage}}%
  \subfigure[Self-avoiding polygon]{
  \begin{minipage}[b]{0.3\textwidth}
    \centering
    \includegraphics{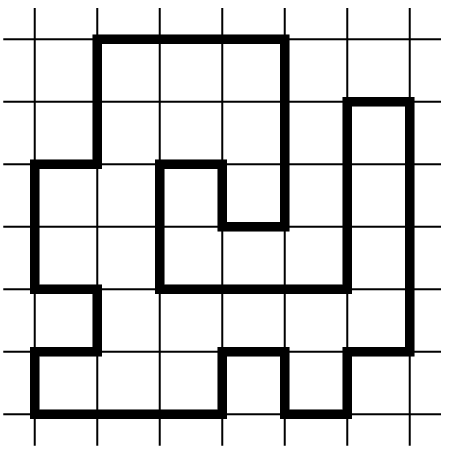}
    \label{subfig:selfavoiding}
  \end{minipage}}
  \caption{Non-convex polyominoes}
\end{figure}

A second class of non-convex polygons is made up of four directed paths.
A {\it three-choice walk} is a self-avoiding walk whose steps are taken in 
accordance with the three-choice rule which allows a step either to the
left or the right or straight ahead after any vertical step, but forbids
a right turn after any horizontal step.  A polygon formed from such
a walk is a {\it three-choice polygon.}  When the walk returns to its starting
point, we don't specify whether the next step, {\em i.e.}, the first step, is 
a valid continuation of the walk.  If it is, the result is a staircase
polygon; if not, it is an {\it imperfect staircase polygon}
(see Figure~\ref{subfig:thchoice}).  When
we refer to three-choice polygons below, we include only
the imperfect ones.

A {\em polyomino} is a union of connected (sharing an edge) cells
of the lattice.  We shall consider one class of nonpolygon polyominoes ---
the {\em staircase polygons with a staircase hole}.  The outer boundary
and the hole are both staircase polygons and must not touch at any point.
An example is shown in Figure~\ref{subfig:spsh}.

\subsection{Self-reciprocity in polyomino enumeration}\label{subsect:examples}
For each of the above classes of  column-convex
polygons, the {\it anisotropic perimeter and area generating function,}
\begin{equation}
G(x,y,q) = \sum_{m\ge1} \sum_{n\ge1} \sum_{a\ge1}
 C(m,n,a) x^m y^n q^a
\label{pagenfun}
\end{equation}
has been computed exactly (see ref.~\cite{bousc} and  references
therein).  Here $C(m,n,a)$ is the number of polygons of the class
with $2m$ horizontal bonds, $2n$ vertical bonds and area $a$.
For the  classes of convex polygons, the anisotropic perimeter generating
function, $G(x,y,1)$ is an algebraic function of the fugacities, $x$ and
$y$, whereas the area generating function, $G(1,1,q)$ is a $q$-series.
For classes of polygons that are only column-convex, both $G(x,y,1)$ and $G(1,1,q)$~\cite{tempa} are
algebraic, but $G(x,y,q)$ involves $q$-series.  A closed-form expression for
the three-choice polygon anisotropic perimeter-area generating function
is not yet known, but by means of a transfer matrix technique it can be
evaluated in polynomial time~\cite{cgd}.  The
isotropic perimeter generating function, $G(x,x,1)$ is known to have
a logarithmic singularity~\cite{cgd}, and is therefore not algebraic, but
is known to be D-finite.  The generating function for staircase
polygons with a staircase hole is also not known in closed form.  Its
properties are expected to be similar in many respects to the generating
function for three-choice polygons~\cite{jge}.

We shall be concerned with self-reciprocity properties of the \gf s
$H_m(y,q)$  that count polygons of width $m$. We
first give two examples.
\begin{enumerate}
\item The area generating function for staircase polygons
of width 4 is the following rational function~\cite{bousb}:
$$H_4(q)= \frac{q^4 (
1+2q+4q^2+6q^3+7q^4+6q^5+4q^6+2q^7+q^8)}{(1-q)^2(1-q^2)^2(1-q^3)^2(1-q^4)}.$$
It satisfies
$$H_4(1/q)= -H_4(q),$$
and is thus self-reciprocal. Observe that the numerator is not only
symmetric (due to self-reciprocity), but also unimodal.
\item The (half-)vertical perimeter and area generating function for
column-convex polygons of width 3 is the  following rational
function, which can be derived from the general formula of ref.~\cite{bousc}:
\begin{align*}
H_3(y,q) =& \frac{yq^3}{(1-yq)^4(1-yq^2)^2(1-yq^3)}\cdot
( y^6q^8+4y^5q^7+2y^5q^6+y^4q^6-y^4q^4 \\
&-4y^3q^5-6y^3q^4-4y^3q^3-y^2q^4+y^2q^2+2yq^2+4yq+1).
\end{align*}
It satisfies
$$H_3(1/y,1/q) = -\frac{1}{yq^3} H_3(y,q)$$
and hence is self-reciprocal. Again, observe the symmetry of the
coefficients in the numerator.
%
%
%
\end{enumerate}
We shall generalize these results to polygons of any width.
Table~\ref{table:summary} summarizes the self-reciprocity properties we have
established. Most of them can be proved in various ways. One can for
instance use a closed form expression of the \gf \ (Section
\ref{subsect:exact}), 
or a functional equation that defines it (Section \ref{sect:temperley});
one can also encode the polygons by a sequence of numbers constrained
by linear diophantine equations and apply Stanley's general results
(Section \ref{sect:stanley}). We shall
see that the last two methods allow us to introduce many additional
parameters and obtain self-reciprocity results that significantly
generalize those of Table~\ref{table:summary}.
\begin{sidewaystable}
\begin{center}
\vskip -5cm
%
\begin{tabular}{|l|c|c|c|}
\hline
Class & Picture & Self Reciprocity & Inversion Relation \\
\hline \hline
&&& \\
Ferrers & \includegraphics[scale=0.2]{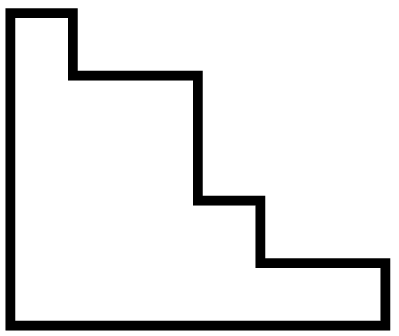} 
& $H_m(1/y,1/q) = (-1)^m y^{m-2} q^{\frac{m^2-3m}{2}} H_m(y,q)$ &
$G(x,y) - y^2G(-x/y,1/y) = 0$ \\
stack & \includegraphics[scale=0.2]{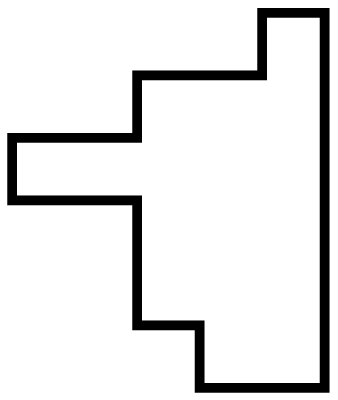} 
& $H_m(1/y,1/q) = -y^{2m-3} q^{m^2-2m} H_m(y,q)$ 
&$G(x,y) + y^3G(x/y^2,1/y) = 0$ 
 \\
staircase &\includegraphics[scale=0.2]{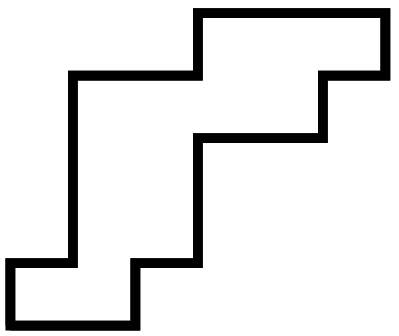}  
& $H_m(1/y,1/q) = -y^{m-1} H_m(y,q), \ m\ge 2$& 
$G(x,y,q) + yG(x/y,1/y,1/q) = -x$ \\
directed convex & \includegraphics[scale=0.2]{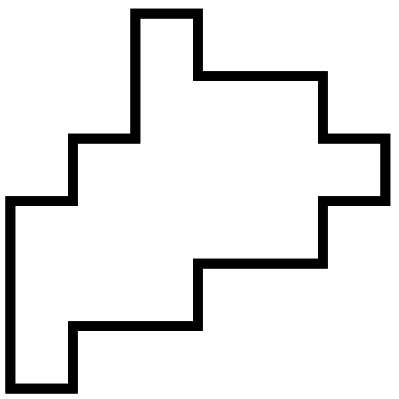} 
%
%
& $H_m(1/y) = -y^{m-2}H_m(y)$
& $G(x,y) + y^2G(x/y,1/y) = 0$ \\
convex & \includegraphics[scale=0.2]{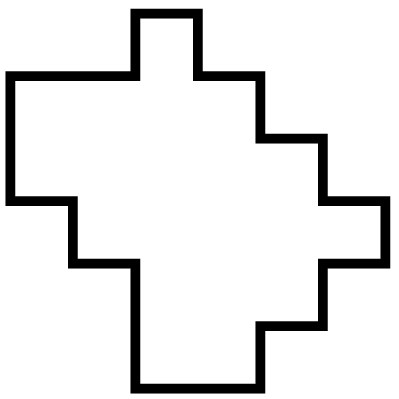} 
& Not simple 
%
%
%
& $G(x,y) + y^3G(x/y,1/y) = xy - x^3 y \frac{\partial}{\partial x}
\frac{1-x+y}{\Delta(x,y)}$\\
%
bargraph & \includegraphics[scale=0.2]{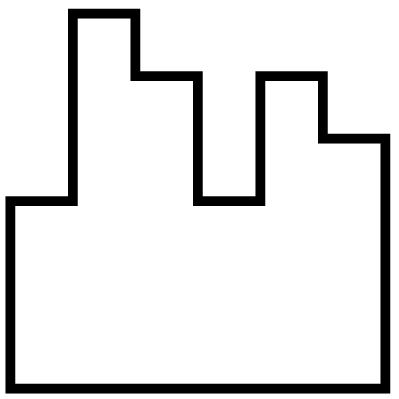} 
& $H_m(1/y,1/q) = \frac{(-1)^m}{yq^m}H_m(y,q)$ &
$ G(x,y,q) - yG(-xq,1/y,1/q) =0$ \\
dir.~col.-conv.
&\includegraphics[scale=0.2]{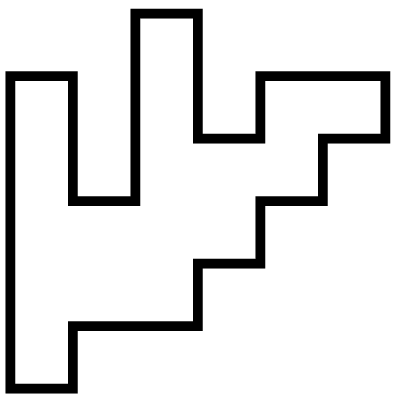} 
%
%
& $H_m(1/q)=-\frac{1}{q}H_m(q)$ & $G(x,q)+qG(x,1/q)=0$ \\
column-convex & \includegraphics[scale=0.2]{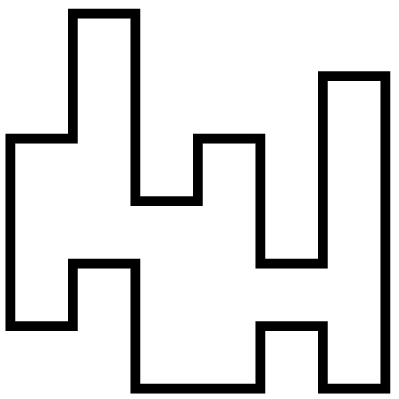} 
& $ H_m(1/y,1/q) = -\frac{1}{y q^m} H_m(y,q) $ &
$G(x,y,q) + yG(xq,1/y,1/q) = 0$\\
&&& \\
\hline
&&& \\
three-choice & \includegraphics[scale=0.2]{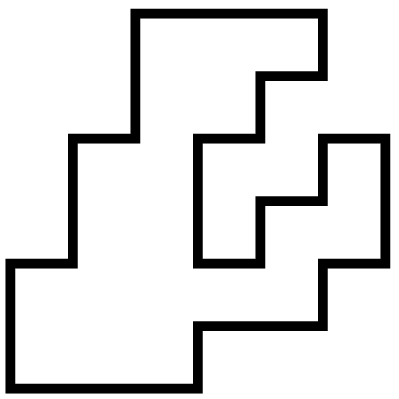}
& Not simple 
& $G(x,y,q) + y^2G(x/y,1/y,1/q) = \mbox{known}$\\
SC with SC hole & \includegraphics[scale=0.2]{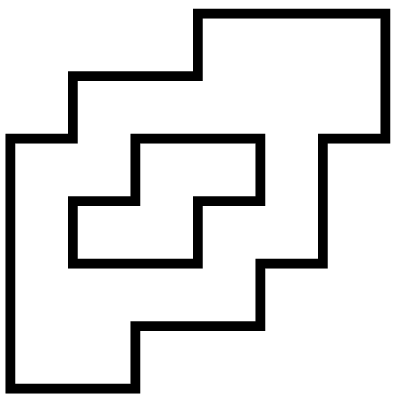} 
& Not simple
& $G(x,y,q) + y^2G(x/y,1/y,1/q) = \mbox{known}$\\
&&& \\
\hline
\end{tabular}
\caption{Summary of polyomino inversion relations}
\label{table:summary}%
\end{center}
\end{sidewaystable}

\subsection{Self-reciprocity via generating functions}\label{subsect:exact}
When a closed form expression for the \gf \ of some class of polygons
is known, it seems  natural to use it to demonstrate an inversion
relation. Let us take the example  of the anisotropic perimeter
generating function for directed convex polygons, which is known to
be~\cite{lin-chang}: 
\begin{equation}
G(x,y) = \frac{xy}{{\sqrt{\Delta(x,y)}}}
\label{dirconclosed}
\end{equation}
with $\Delta(x,y) =
1-2x-2y-2xy+x^2+y^2=(1-y)^2\left[1-x(2+2y-x)/(1-y)^2\right]$. 
Expanding this expression in $x$ gives
\begin{eqnarray*}
G(x,y) &=& \sum_{m\ge1} H_m(y) x^m \\
&=& \frac{y}{1-y} x +
\frac{y(1+y)}{(1-y)^3} x^2 + \frac{y(1+4y+y^2)}{(1-y)^5}x^3 +
\frac{y(1+9y+9y^2+y^3)}{(1-y)^7}x^4+ O(x^5)
\end{eqnarray*}
which suggests that the partial generating functions, $H_m(y)$ are
self-reciprocal, and more precisely, that 
%
%
$H_m(1/y) = -y^{m-2}H_m(y)$. This is equivalent to the inversion relation 
\begin{equation}
G(x,y)+ y^2 G(x/y,1/y) = 0,
\label{dirconinv}
\end{equation}
which is easily checked from the
closed form of the generating function.
Note that an explicit expression for $H_m(y)$ is given in~\cite{bousb}.
%
The inversion relations for convex polygons 
and directed column-convex
polygons may also be obtained from
the expression of their \gf .

The partial generating functions for directed convex polygons,
counted by the area, are
not  self-reciprocal: for instance, the \gf \ for width $3$ is
$$q^3(1+3q+3q^2+2q^3+q^4)/(1-q)^2(1-q^2)^2(1-q^3).$$
However, many other classes of column-convex polygons have an inversion
relation for the full anisotropic perimeter and area generating function.
Since these generating functions are also known in closed form they could
be derived as above.  However more can be shown, namely that there is a
self-reciprocity for any parameter which is a linear function of the
{\em vertical heights\/} in the graph.  This very general result will be derived in
Section~\ref{sect:temperley}.
Likewise, the inversion
relations for three-choice polygons and staircase polygons with a
staircase hole, given in
Table~\ref{table:summary}, are also
special cases of more general formulae which will be derived in
Section~\ref{sect:stanley}.

\subsection{Using inversion relations to compute\\ generating functions}
\label{subsect:computegf}
As in statistical mechanics, the inversion relation and symmetry, and some
general assumptions on analyticity of the generating function,
are sometimes sufficient
to determine the solution completely.
In order to have an algorithm for computing a generating function term
by term, it is necessary, but not sufficient,
to have some property relating terms of different
orders.  For our purposes this property will always be $x$-$y$ symmetry.
Thus we restrict our attention  to
classes of graphs with $x$-$y$ symmetry, {\em i.e.}, 
Ferrers, staircase, directed
convex, convex and three-choice polygons, and staircase polygons with
a staircase hole.  Moreover, we shall only consider  the anisotropic
perimeter generating 
function (without area). For the former four classes we will show that
the inversion relation provides sufficient additional information to
compute the generating function, whereas for the latter two it does
not.

The general form of the generating function is
\begin{equation}
G(x,y) =  H_1(y) x + H_2(y) x^2 + H_3(y) x^3+\cdots
\label{Gexp}
\end{equation}
where the partial generating functions, $H_m(y)$ are rational functions
\begin{equation}
H_m(y)=\frac{P_m(y)}{D_m(y)},
\label{Hexp}
\end{equation}
with $D_m(0)=1$. The general form of the inversion relation is
\begin{equation}
G(x,y) \pm y^\alpha G(\epsilon x/y,1/y) = \text{RHS}
\label{generalinv}
\end{equation}
where $\alpha$ is an integer and RHS is zero or some simple function.
It is equivalent to a self-reciprocity relation of the form
$$ H_m(1/y)\pm \epsilon^m y^{m-\alpha}H_m(y)=\text{RHS}.$$
Whether the inversion relation is sufficient to compute the generating
function depends on the value of the exponent $\alpha$
and on the degree of the denominator, $D_m(y)$.
Direct proof of the denominator form
can often be obtained.  
For Ferrers graphs, it is easily shown that $D_m(y)=(1-y)^{m}$. For
staircase polygons, one finds
$D_m(y)=(1-y)^{2m-1}$.  The same denominator form holds for
directed convex and convex polygons also.
For the three-choice polygons and staircase
polygons with a staircase hole it can
be shown that the denominators are
\begin{equation}
D_m(y) = \begin{cases} (1-y)^{2m-1} (1+y)^{2m-7} & \text{$m$ even} \\
(1-y)^{2m-1} (1+y)^{2m-8} & \text{$m$ odd.} \end{cases}
\label{threechoiceden}
\end{equation}
We assume that in general we know the denominator form either empirically or by
rigorous proof, and that $D_m(y)$ is of degree $d_m$.

Now we proceed inductively, following Baxter~\cite{baxb}.  If we have
already computed the coefficient functions $H_1(y),\ldots,H_{m-1}(y)$ in
the expansion~(\ref{Gexp}) and if $x$-$y$ symmetry holds, we also
know the coefficients of $y,y^2,\ldots,y^{m-1}$ in the expansion of $G(x,y)$.
In particular, we can compute the coefficients of $y,y^2,\ldots,y^{m-1}$
in the numerator polynomial $P_{m}(y)$.  In order to obtain the unknown
coefficients of $P_{m}(y)$, we must be able to express them in terms
of the known ones by means of the inversion relation.  Writing
$P_m(y)=\sum_k a_k y^k$, and using $D_m(1/y)=\pm y^{-d_m}D_m(y)$, the
inversion relation fixes the value of the combinations of coefficients,
$a_k\pm a_\ell$, with $k+\ell=\alpha+d_m-m$. Hence the
determination of all the coefficients $a_k$ is possible
if and only if the arithmetic condition holds:
\begin{equation}
d_{m} < 3m - \alpha.
\label{arthcond}
\end{equation}
This condition is seen to hold for
all the classes of convex polygons we have looked at, since $d_m\le 2m-1$,
but not for three-choice polygons or staircase polygons with a staircase
hole, since $d_m\sim 4m$. 
%
%

\section{Self-reciprocity via Temperley methodology}
\label{sect:temperley}
\begin{figure}
  \centering
  \psfrag{gt}{$\overline\gamma$}
  \psfrag{gb}{$\underline\gamma$}
  \includegraphics{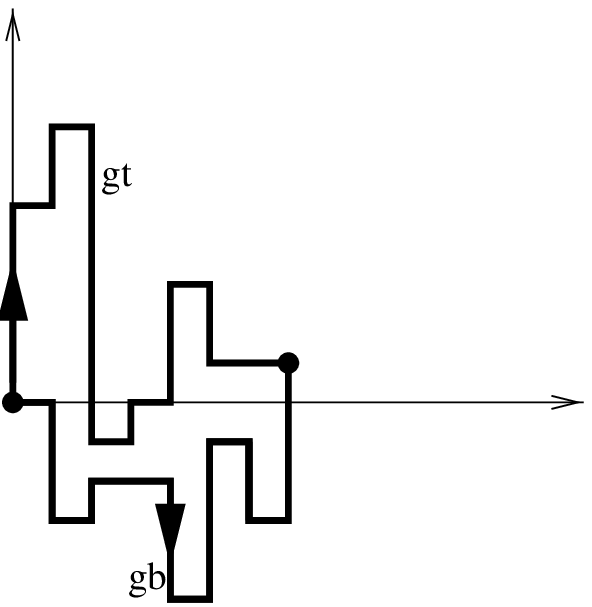}
  \caption{A column-convex polygon}
\label{fig:vc}
\end{figure}

We consider column-convex polygons as pairs of partially directed paths
having the same endpoints, as indicated in Figure \ref{fig:vc}.

Let $P$ be a column-convex polygon of width $m$. 
For $0 \le i \le m$,  we denote by  $\Nt_i$ (resp. $\St_i$) the
number of north  (resp. south)
steps in the top path $\overline\gamma$ at abscissa $i$. 
For $0 \le i \le m$, we denote by  $\Nb_i$
(resp. $\Sb_i$) the number of north (resp. south)
steps of the bottom path $\underline\gamma$ at abscissa $i$.
We choose the end points of the paths in such a way that
$\Nb_0=\Sb_0=\Nt_m=\St_m=\St_0=\Sb_m=0$.
Note that 
$$
\sum_{k=0}^{m}(\Nt_k+\Sb_k-\St_k-\Nb_k)= 0.$$

We notice that  all standard statistics are {\em linear\/}
functions of the $ \Nt_i, \St_i, \Nb_i$ and $\Sb_i$. For instance,
the vertical perimeter of the polygon is
\begin{align}
2n&= \sum_{k=0}^{m}(\Nt_k+\Sb_k+\St_k+\Nb_k)\notag\\
&= 2\, \sum_{k=0}^{m}( \Nt_k+\Sb_k).
\label{vperimeter}
\end{align}
The height of the $i^{\text{th}}$ column of the
polygon is, for $ 1 \le i \le m$,
\begin{align*}
h_i & 
 = \sum_{k=0}^{i-1} (\Nt_k+\Sb_k-\St_k-\Nb_k),
\end{align*}
and the area of the polygon is
\begin{equation}
a= \sum_{k=0}^{m} (m-k)(\Nt_k+\Sb_k-\St_k-\Nb_k).
\label{area}
\end{equation}

\begin{theo} \label{theo:reci}
Let $\cal P$ be one of the following sets:
 Ferrers diagrams,
 stacks (drawn as in Figure~\ref{subfig:stackh}),
 staircase polygons,
 bar-graphs,
 column-convex polygons. 
Let ${\cal P}_m$ be the subset
of $\cal P$ containing all polygons of
width $m$. Let $F_m$ be the \gf \ for polygons in the set ${\cal P}_m$:
$$F_m(\byt,\bzt,\byb,\bzb) = \sum _{P \in {\cal P}_m}
\byt^\bNt \bzt^\bSt \byb^\bNb \bzb^\bSb. $$
Then $F_m$ is a rational function, and it is self-reciprocal:
\begin{equation}
F_m(1/\byt, 1/\bzt, 1/\byb, 1/\bzb) =
 C_m F_m(\byt, \bzt, \byb, \bzb),
\label{ccreci}
\end{equation}
with 
$$C_m = \left\{ \begin{array}{ll}
\displaystyle \frac{(-1)^{m}\yb_m^{m-2} }{\yt_0}\prod_{i=1}^{m-1} \yt_i  
& \hbox{for Ferrers graphs,} \\
& \\
\displaystyle
- \frac{\yb_m^{2m-3}}{\yt_0} \prod_{i=1}^{m-1} \yt_i  \zb_i  & \hbox{for stacks,}
\\
& \\
\displaystyle -\prod_{i=1}^{m-1} \yb_i  \yt_i & \hbox{for staircase polygons }
 (m\ge2),\\
& \\
\displaystyle \frac{(-1)^{n}}{\yt_0\yb_m}      & \hbox{for bar-graphs},\\
& \\
\displaystyle
-\frac{1}{\yt_0 \yb_m} & \hbox{for column-convex polygons}. \\
\end{array} \right.$$
\end{theo}

The proof of the theorem is based on the so-called Temperley approach
for counting column-convex polygons~\cite{tempa}, combined with
the systematic use of formal power series~\cite{bousc}. Here we
provide only the proof for column-convex polygons, since the others
are very similar.

We commence by showing that the partial generating functions for
column-convex polygons, $V_m(\byt, \bzt, \byb, \bzb)$, can be computed
inductively.

\begin{propo}\label{functional}
Let $V_m(\byt, \bzt, \byb, \bzb)$ be
the \gf \ for column-convex
polygons of width $m$. Let us denote it, for the sake of simplicity,
$V_m(\yb_m)$.   Then the series $V_m(\yb_m)$ can be defined inductively by:
$$V_1(\yb_1)=\frac{\yt_0\yb_1}{1-\yt_0\yb_1}$$
and
\begin{align*}
&V_{m+1}(\yb_{m+1}) = \frac{(1-\yb_m\zb_m)(1-\yt_m\zt_m)V_{m}(\yb_{m+1})}
{(1-\yb_{m+1}\yt_m)(1-\yb_{m+1}\zb_m)(1-\yb_{m+1}^{-1}\yb_m)(1-\yb_{m+1}^{-1}\zt_m)} \\
&\quad + \frac{(\zt_m -\yb_{m+1} \yb_m \zb_m)
V_m(\zt_m)}{(1-\yb_{m+1}\zb_m)(1-\yb_{m+1}^{-1}\zt_m)(\yb_m-\zt_m)}
+ \frac{(\yb_m -\yb_{m+1} \yt_m \zt_m)
V_m(\yb_m)}{(1-\yb_{m+1}\yt_m)(1-\yb_{m+1}^{-1}\yb_m)(\zt_m-\yb_m)}.
\end{align*}
\end{propo}

\noindent
{\bf Proof.} The basic idea is build a polygon of width $m+1$ by
adding a new column to a polygon of width $m$ \cite{bousc}.  It is
convenient to use Hadamard products to establish the functional
equation.

Let $F(t)=\sum f_ht^h$ and $R(t)=\sum r_ht^h$ be two formal power
series in $t$ with coefficients in a ring $A$. We denote by $F(t)\odot R(t)$
the Hadamard product of $F(t)$ and $R(t)$, evaluated at $t=1$:
$$F(t) \odot R(t) = \sum f_hr_h.$$
In what follows, $f_h$ (resp.\ $r_h$) will be the generating function
for some column-convex polygons whose rightmost (resp.\ leftmost) column has
height $h$, so that $F(t) \odot R(t)$ will count polygons obtained by
{\em matching\/} the rightmost column of a polygon of type $F$ with the
leftmost column of a polygon of type $R$. Also, $R(t)$ will be a
rational function of $t$. We shall 
use the following simple identity:
$$F(t)\odot \frac{1}{1-at}=F(a).$$

\medskip

The expression for $V_1(\yb_1)$ is obvious. We build a column-convex
polygon of width $m+1$ as follows: we take a polygon of width $m$ and
match its rightmost column with the leftmost column of
a column-convex polygon of width $2$. 

\begin{figure}
  \centering
  \psfrag{gt0}{$>0$}
  \psfrag{ge0}{$\ge0$}
  \psfrag{$+$}{$+$}
  \psfrag{O}{$\odot$}
  \psfrag{$m$}{$h$}
  \includegraphics{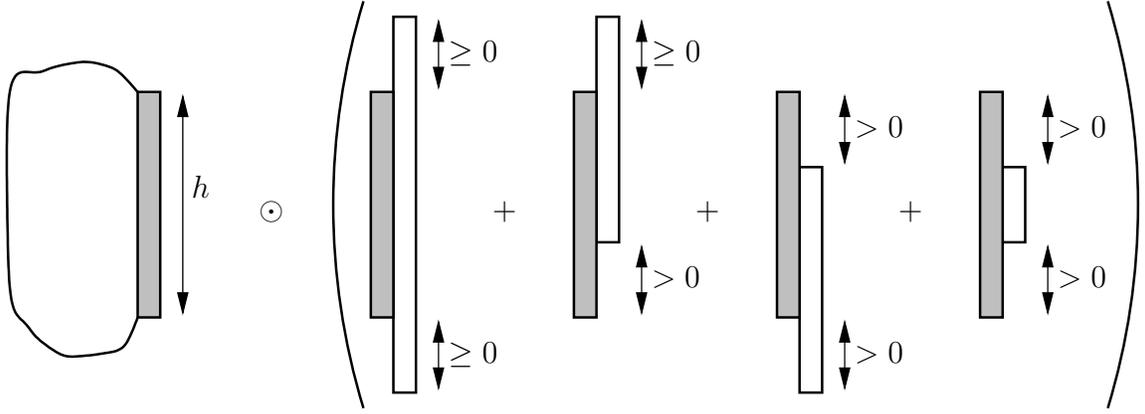}
  \caption{Construction of column-convex polygon by Hadamard products.}
\label{fig:hadamard}
\end{figure}
This is illustrated by Figure~\ref{fig:hadamard}, which shows that
\begin{equation}
 V_{m+1}(\yb_{m+1})= V_m(t)\odot R(t),
\label{rec}
\end{equation}
where
\begin{align*}
R(t) =& \frac{t\yb_{m+1}}{1-t\yb_{m+1}} \cdot \frac{1}{1-\yb_{m+1}\yt_m}
 \cdot \frac{1}{1-\yb_{m+1}\zb_m}
+ \frac{t\yb_{m+1}}{1-t\yb_{m+1}} \cdot \frac{1}{1-\yb_{m+1}\yt_m}
 \cdot \frac{t\yb_{m}}{1-t\yb_m} \\
&+ \frac{t\yb_{m+1}}{1-t\yb_{m+1}} \cdot  \frac{t\zt_m}{1-t\zt_m}
 \cdot \frac{1}{1-\yb_{m+1}\zb_m}
+ \frac{t\yb_{m+1}}{1-t\yb_{m+1}} \cdot \frac{t\zt_m}{1-t\zt_m}
 \cdot \frac{t\yb_m}{1-t\yb_m}
\end{align*}
is the \gf \ for column-convex polygons of width $2$.
In order to determine the coefficient $r_h$ of $t^h$ in $R(t)$, we
expand $R(t)$ in partial fractions of $t$:
\begin{align*}
R(t)=& \frac{\zb_m \yt_m \yb_{m+1}^2}{(1-\yb_{m+1}\yt_m)(1-\yb_{m+1}\zb_m)} \\
&+
\frac{(1-\yb_m \zb_m)(1-\yt_m \zt_m)}
{(1-\yb_{m+1} \yt_m)(1-\yb_{m+1}\zb_m)(1-\yb_{m+1}^{-1}\yb_m)
 (1-\yb_{m+1}^{-1} \zt_m)} \cdot \frac{1}{1-t\yb_{m+1}} \\
&+
 \frac{(\zt_m-\yb_{m+1}\yb_m \zb_m)}{(1-\yb_{m+1} \zb_m)(1-\yb_{m+1}^{-1}\zt_m)
 (\yb_m-\zt_m)} \cdot \frac{1}{1-t \zt_m} \\
&+
 \frac{(\yb_m-\yb_{m+1} \yt_m \zt_m)}{(1-\yb_{m+1} \yt_m)(1-\yb_{m+1}^{-1}\yb_m)
 (\zt_m-\yb_m)} \cdot \frac{1}{1-t \yb_m}.
\end{align*}
Note that $V_m(0)=0$.
We now combine eqn.~\Ref{rec} with the above  expression for $R(t)$ to
obtain the announced expression for $V_{m+1}(\yb_{m+1})$.
\qed

\bigskip

\noindent {\bf Proof of Theorem \ref{theo:reci}.} 
 Induction on $m$ using the functional equation of
Proposition~\ref{functional} shows that the partial generating
functions for column-convex polygons satisfy:
$$ {V_m}(1/\byt, 1/\bzt, 1/\byb,1/\bzb) = -\frac{1}{\yt_0 \yb_m} V_m(\byt, \bzt, \byb,\bzb).$$
\qed

We proceed similarly for the other families: the functional equation
is obtained by setting some of the variables $\yt_i, \yb_i, \zt_i$ and
$
\zb_i$ to $0$.  Then, an inductive argument yields the
self-reciprocity result. 

It would be tempting to write that the self-reciprocity of $V_m$
implies the self-reciprocity of, say, the \gf \ for staircase
polygons, obtained by setting $\zb_i$ and $\zt_i$ to $0$ in $V_m$. But
replacing a variable by $0$ in a self-reciprocal
rational function might break the self-reciprocity: for instance, take
$P(y_1,y_2)=1+y_1+2y_1^2+2y_2+y_1y_2+y_1^2y_2.$ Then
$$P(1/y_1,1/y_2)=\frac{1}{y_1^2y_2}P(y_1,y_2),$$
but $P(y_1,0)=1+y_1+2y_1^2$ is {\em not\/} self-reciprocal.

\medskip

However, the following simple lemma gives a useful stability property of
self-reciprocal rational functions.
\begin{lem} \label{lemma:stability}
Let $F(y_1,\ldots,y_n)$ be a self-reciprocal rational function.
Let $\ssA$ be an $m\times n$ integer matrix.  Let
$\bu=(u_1,\ldots,u_m)$, and define $\bu^\ssA$ to be the $n$-tuple
whose $i^{\text{th}}$ coordinate is
$\prod_k u_k^{a_{ki}}$.  Then the series $G(\bu)=F(\bu^\ssA)$, if defined,
is self-reciprocal in the variables $u_i$.  More precisely,
if $F(1/\by)=\pm \by^\bbeta F(\by)$ then
$G(1/\bu)=\pm \bu^{\ssA\bbeta} G(\bu)$.
\end{lem}

From Theorem~\ref{theo:reci} and Lemma~\ref{lemma:stability} we immediately
deduce:
\begin{coro} \label{coro:linear}
For any of the sets $\cal P$ listed in Theorem~\ref{theo:reci}, and any
statistics on column-convex polygons that can be expressed as {\em
linear} functions of the quantities $\bNt,\bSt,\bNb,\bSb$,
the \gf \ for
polygons in the set ${\cal P}_m$ according to these statistics is a
self-reciprocal rational function.
\end{coro}

This corollary allows us to complete the top part of
Table~\ref{table:summary}.  Let
us, for instance, derive the inversion relation satisfied by the
tri-variate \gf \ $G(x,y,q)$ for  column-convex polygons,  taking into account
the usual parameters of interest: horizontal and vertical half-perimeters (variables $x$ and $y$), and area (variable $q$).

Eqns.~\Ref{vperimeter} and \Ref{area} express the vertical perimeter and the area
in terms of the quantities $\bNt, \bNb, \bSt$ and $\bSb$. They imply
that the (half) vertical perimeter and area \gf \ $H_m(y,q)$ for
column-convex polygons of width $m$ is
$$H_m(y,q)= V_m(\bNt, \bSt, \bNb,  \bSb)$$
where 
 $\yt_k = \zb_k= yq^{m-k}$ and $ \yb_k=\zt_k= q^{-(m-k)}$.
Theorem \ref{theo:reci} then gives
$$H_m(1/y,1/q)=-\frac{1}{yq^m}H_m(y,q),$$
which implies
\begin{equation}
G(x,y,q)+y G(x q,1/y,1/q) = 0.
\label{colconvinv}
\end{equation} 

Note that in the first two self-reciprocity relations of Table 1, the
exponent of $q$ depends quadratically on the width.  For this reason,
they only yield an inversion relation for $q=1$.

%
\section{Self-reciprocity via Stanley's general results}
\label{sect:stanley}
\subsection{Linear homogeneous diophantine systems}\label{subsect:lhds}
Stanley has analyzed the situation where the objects
to be counted correspond to integer solutions of a system of
linear equations with integer coefficients (linear diophantine system)
subject to
a set of constraints.  He has established certain
conditions under which reciprocity
relations will hold between two combinatorics problems defined by the same
linear diophantine system but by different sets of constraints, and also
conditions
under which the solution to a given problem will be
self-reciprocal~\cite{stanb,stane}.

Consider the linear homogeneous diophantine system (LHD-system),
\begin{equation}
\bPhi\balpha=\bzero
\label{lhds}
\end{equation}
in the unknowns $\balpha=(\alpha_1,\ldots,\alpha_s)$ with $\bPhi$ a matrix of integers having
$p$ rows and $s$ columns and $\bzero$ a $p$-tuple of zeros.
The corank $\kappa$
of the system is defined to be $s-\mathrm{rank}(\bPhi)$.
For a linearly independent system, $\kappa=s-p$.  Let $S$ be a set of
integer solutions to eqn.~(\ref{lhds}).  We define the generating function,
$S(\by)$, as the formal power series
\begin{equation}
S(\by) = \sum_{\balpha\in S} \by^\balpha
\label{lhdgf}
\end{equation}
where $\by=(y_1,\ldots,y_s)$ is a vector of fugacities associated with
the unknowns in eqn.~(\ref{lhds}).

In our applications, we find two types of constraints on the unknowns, $\alpha_j$.
Certain of the unknowns, $\alpha_j$, are
required to be strictly positive while the rest
are required to be non-negative.
Conveniently,
precisely these kinds of constraints have been treated by Stanley.  
Let the unknowns be
%
%
$\balpha=(\bgamma,\bdelta)\equiv\bgamma\oplus\bdelta$
where $\bgamma$ is an $n$-tuple and $\bdelta$ is an $(s-n)$-tuple.  Likewise
let $\by=(\bu,\bv)$.  In what follows, the notation, $\bdelta>\bzero$, means
that all coordinates of $\bdelta$ are positive.
\begin{propo}
Let $E$ be the set of integer
solutions, $(\bgamma,\bdelta)$, to a linear
homogeneous diophantine system of corank $\kappa$, such that
$\bgamma\ge\bzero$ and $\bdelta>\bzero$.
Let $\overline E$ be the set of solutions to the same system with
$\bgamma>\bzero$ and $\bdelta\ge\bzero$.  If the system has an integer
solution, $(\bgamma,\bdelta)$ such that $\bgamma>0$ and
$\bdelta<0$, then $E(\bu,\bv)$ and $\overline E(\bu,\bv)$ are rational
functions obeying the reciprocity relation
\begin{equation}
\overline E(\bu,\bv) = (-1)^\kappa E(1/\bu,1/\bv).
\label{lhdrecb}
\end{equation}
\label{propo:recb}
\end{propo}
{\em Proof.}  This is Proposition 8.3 of ref.~\cite{stanb} and the proof
is given there. \qed

Proposition~\ref{propo:recb} can be specialized
to obtain a self-reciprocity condition,
which will be our main tool in the derivations to follow.
\begin{coro}
A sufficient condition for the function $E(\bu,\bv)$ to be self-reciprocal
is that the linear homogeneous diophantine system has the solution
$(\bgamma,\bdelta)=(\bone,-\bone)$.  In this case
\begin{equation}
E(1/\bu,1/\bv)=(-1)^\kappa \frac{\bu^\bone}{\bv^\bone} E(\bu,\bv).
\label{lhdselfb}
\end{equation}
\label{coro:selfb}
\end{coro}
{\em Proof.} Since the solution $(\bone,-\bone)$ satisfies the conditions
of Proposition~\ref{propo:recb}, the reciprocity result~(\ref{lhdrecb}) holds.
The result follows immediately
from the shift 
$(\bgamma,\bdelta)\rightarrow(\bgamma+\bone,\bdelta-\bone)$ which establishes a bijection
between the sets $E$ and $\overline E$. \qed

Since the conditions of the corollary are sufficient but not necessary,
it is often possible to find a perfectly valid 
LHD-system  describing a given self-reciprocal 
generating function, $E(\bu,\bv)$, which does {\em not} admit the solution
$(\bone,-\bone)$.
Hence we are faced with the problem
of finding a suitable LHD-system which satisfies the corollary.  A
useful heuristic is to start with
an LHD-system in many unknowns, and selectively eliminate those unknowns
whose constraints are not independent of the constraints on the other
unknowns.  In all the cases we will consider,
the resulting system will satisfy the conditions of Corollary~\ref{coro:selfb}.
We do not justify  this heuristic here.
In a paper subsequent to ref.~\cite{stanb},
Stanley~\cite{stand} develops a more comprehensive theory which overcomes
these difficulties, and which additionally gives ``correction'' terms
for systems in which self-reciprocity fails to hold.  We have not yet
explored the ramifications of this theory.

Before applying the above result to staircase polygons with a staircase hole or
to three-choice polygons, we use it to derive the reciprocity relation for
ordinary staircase polygons of width three.  This will serve to
illustrate all the
basic ingredients of the method.

\bigskip
\noindent
\begin{figure}
  \centering
  \psfrag{R1}{\small $\Nt_1$}
  \psfrag{R2}{\small $\Nt_2$}
  \psfrag{S0}{\small $M_0$}
  \psfrag{S1}{\small $M_1$}
  \psfrag{S2}{\small $M_2$}
  \psfrag{S3}{\small $M_3$}
  \psfrag{T1}{\small $\Nb_1$}
  \psfrag{T2}{\small $\Nb_2$}
  \includegraphics{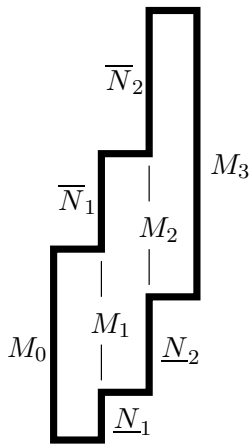}
  \caption{Staircase polygon of width three}
\label{fig:stairthree}
\end{figure}
\begin{exam}
Staircase polygons of width three can be characterized
by the heights  $\Nt_1$,
$\Nt_2$, $\Nb_1$, $\Nb_2$,  $M_0$, $M_1$, $M_2$ and  $M_3$, as shown in Figure
\ref{fig:stairthree}.
Decomposing the polygon
into three columns, and imposing the condition that each
column be as high on the left as it is on the right, we obtain
the linear homogeneous diophantine system
\begin{subequations}
\begin{align}
M_0-M_1-\Nb_1 &= 0 \label{exlhda} \\
M_1+\Nt_1-M_2-\Nb_2 &= 0 \label{exlhdb} \\
M_2+\Nt_2-M_3 &= 0. \label{exlhdc}
\end{align}
\end{subequations}
All heights must be nonnegative, but the self-avoidance condition additionally
requires that the $M_j$ be positive.  The constraints $M_0>0$ and $M_3>0$
are actually redundant, since they follow from
eqns.~(\ref{exlhda},\ref{exlhdc}) and the
constraints on the remaining 
unknowns, namely
\begin{align}
\Nt_1,\Nt_2,\Nb_1,\Nb_2 &\ge 0 \notag\\
M_1,M_2 &> 0.
\label{exconst}
\end{align}
Since the constraints on $M_0$ and $M_3$ play no role in the solution,
we are free to eliminate these unknowns, and it turns out to be necessary
to do so in order to apply Corollary~\ref{coro:selfb}.  We are left with
the single equation~(\ref{exlhdb})
in the six independent unknowns $\bgamma=(\Nt_1,\Nt_2,\Nb_1,\Nb_2)$
and $\bdelta=(M_1,M_2)$.
Let us associate to the unknown $\Nt_i$ (resp. $\Nb_i$, $M_i$) the
fugacity $\yt_i$ (resp. $\yb_i$, $z_i$).

Let $E'$ be the set of solutions to eqn.~(\ref{exlhdb}) subject to the
the constraints $\bgamma\ge0$ and $\bdelta>0$.  
Since $\bgamma=\bone$, $\bdelta=-\bone$ is a solution to eqn.~(\ref{exlhdb}),
Corollary~\ref{coro:selfb} tells us that $E'(\byt,\byb,\bz )$ is self-reciprocal,
\begin{equation}
E'(1/\byt,1/\byb,1/\bz ) = -\frac{\yt_1 \yt_2 \yb_1 \yb_2}{z_1 z_2}
E'(\byt,\byb,\bz ).
\label{exselfa}
\end{equation}

Equations~(\ref{exlhda},\ref{exlhdc}) imply that to account for
the dependent parameters $M_0$ and $M_3$, we make the
substitutions $z_1\rightarrow z_0 z_1$, $\yb_1\rightarrow z_0 \yb_1$,
$z_2\rightarrow z_2 z_3$ and $\yt_2\rightarrow z_3 \yt_2$.
Applying Lemma~\ref{lemma:stability}, we obtain for the set $E$ of
%
%
nonnegative solutions to ~(\ref{exlhda},\ref{exlhdb},\ref{exlhdc})
such that $\bM >{\bzero}$:
\begin{equation}
E(1/\byt,1/\byb,1/\bz ) = -\frac{\yt_1 \yt_2 \yb_1 \yb_2}{z_1 z_2}
E(\byt,\byb,\bz ).
\label{exselfb}
\end{equation}
Notice that reintroducing the dependent unknowns has not changed
the constant factor.  This feature holds as well in the more complicated
models we will look at.
The result~(\ref{exselfb}) may be verified by inspection of the explicit
expression for the generating function
$$
E(\byt,\byb,\bz)=
 \frac{z_0 z_1 z_2 z_3 (1-\yt_1 z_0 z_1 z_2 z_3 \yb_2)}
 {(1-z_0 \yb_1) (1-z_0 z_1 \yb_2) (1-\yt_1 \yb_2)(1-z_0 z_1 z_2 z_3)
 (1-\yt_1 z_2 z_3)(1-\yt_2 z_3)}.
$$
\end{exam}

\subsection{Applications}\label{subsect:applications}
We now apply the methods of Section~\ref{subsect:lhds} to staircase polygons
with a staircase hole and to three-choice polygons.  All
the essential steps have already been
seen in the derivation of the reciprocity result for staircase polygons
of width three.  They are
\begin{enumerate}
\item Set up a linear homogeneous diophantine system by decomposing the
polyomino into width one rectangles and imposing the condition that the
left and right sides of each rectangle have equal height.
\item Sort the unknowns into three classes, $\bgamma$, $\bdelta$ and $\btau$,
according to whether they are constrained to be nonnegative, constrained
to be positive or constrained by conditions on the other unknowns.
\item Use Gaussian elimination to remove the unknowns in $\btau$.
\item Verify that the resulting system is solved by setting all members of
$\bgamma$ equal to one and all members of $\bdelta$ equal to minus one.  Apply
Corollary~\ref{coro:selfb} to obtain the self-reciprocity result for the
reduced system.
\item Reintroduce the unknowns in the set $\btau$ by
means of Lemma~\ref{lemma:stability}.
\end{enumerate}

\begin{figure}
  \centering
  \psfrag{R1}{\small $\Nt_1$}
  \psfrag{R2}{\small $\Nt_2$}
  \psfrag{R3}{\small $\Nt_3$}
  \psfrag{R4}{\small $\Nt_4$}
  \psfrag{R5}{\small $\Nt_5$}
  \psfrag{R6}{\small $\Nt_6$}
  \psfrag{S0}{\small $M_0$}
  \psfrag{S1}{\small $M_1$}
  \psfrag{S2}{\small $M_2$}
  \psfrag{S6}{\small $M_6$}
  \psfrag{S7}{\small $M_7$}
  \psfrag{T1}{\small $\Nb_1$}
  \psfrag{T2}{\small $\Nb_2$}
  \psfrag{T3}{\small $\Nb_3$}
  \psfrag{T4}{\small $\Nb_4$}
  \psfrag{T5}{\small $\Nb_5$}
  \psfrag{T6}{\small $\Nb_6$}
  \psfrag{U3}{\small $\Ht_3$}
  \psfrag{U4}{\small $\Ht_4$}
  \psfrag{U5}{\small $\Ht_5$}
  \psfrag{V2}{\small $M_2$}
  \psfrag{V3}{\small $M_3$}
  \psfrag{V4}{\small $M_4$}
  \psfrag{V5}{\small $M_5$}
  \psfrag{W3}{\small $\Hb_3$}
  \psfrag{W4}{\small $\Hb_4$}
  \psfrag{X2}{\small $\Mt_2$}
  \psfrag{X3}{\small $\Mt_3$}
  \psfrag{X4}{\small $\Mt_4$}
  \psfrag{X5}{\small $\Mt_5$}
  \psfrag{X6}{\small $\Mt_6$}
  \psfrag{Z2}{\small $\Mb_2$}
  \psfrag{Z3}{\small $\Mb_3$}
  \psfrag{Z4}{\small $\Mb_4$}
  \psfrag{Z5}{\small $\Mb_5$}
  \subfigure[Staircase polygon with staircase hole, $(k,\ell,m)=(3,5,7)$]{
  \begin{minipage}[b]{0.45\textwidth}
    \centering
    \includegraphics{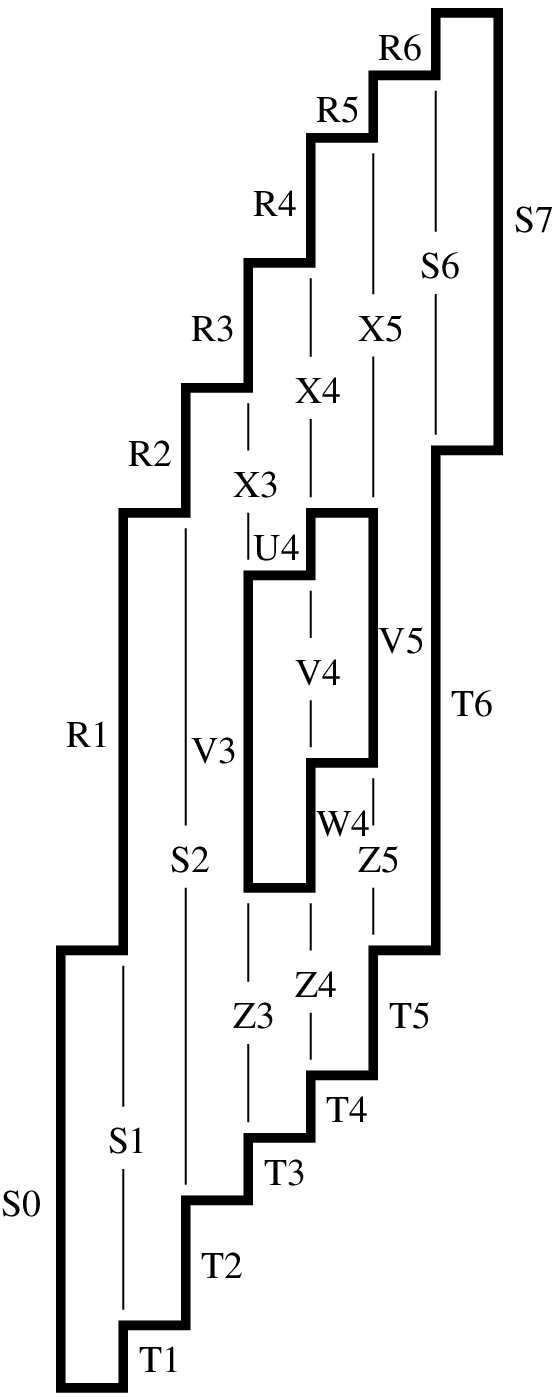}
    \label{subfig:spshb}
  \end{minipage}}%
  \subfigure[Three-choice polygon, $(k,\ell,m)=(2,6,5)$]{
  \begin{minipage}[b]{0.45\textwidth}
    \centering
    \includegraphics{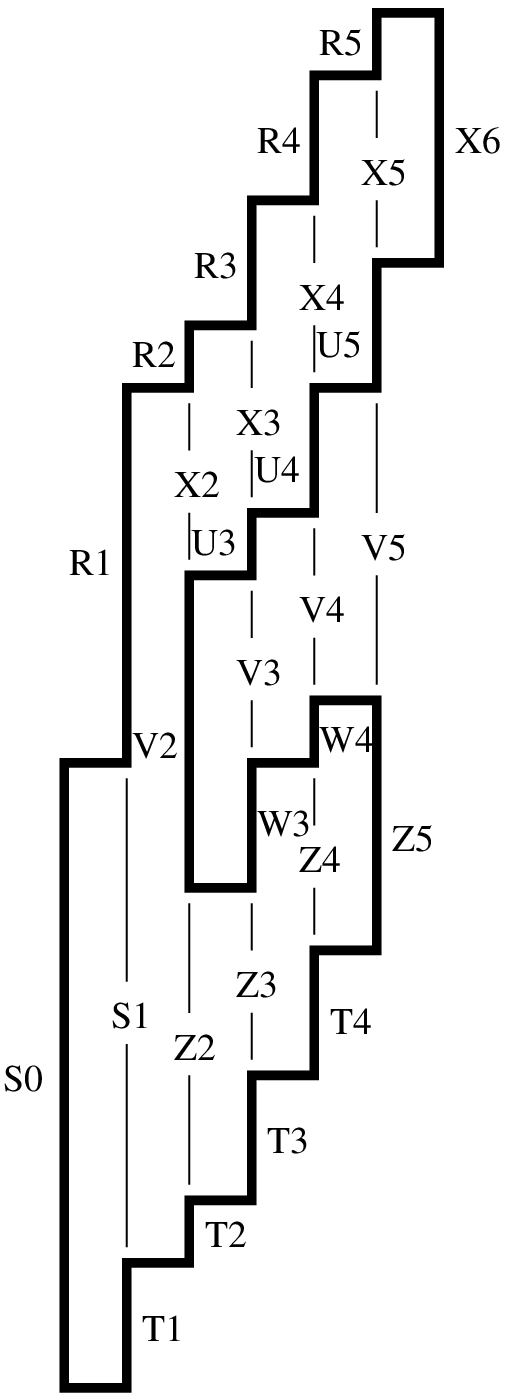}
    \label{subfig:threechoice}
  \end{minipage}}
  \caption{Labels for polyomino vertical heights}
  \label{fig:heights}
\end{figure}
We can define three widths for a staircase polygon with a staircase
hole:
the distance from the left edge of the figure to the
left edge of the hole, $k$, the distance from the left edge of the figure
to the right edge of the hole, $\ell$, and the width of the entire figure, $m$.
Note that $0<k<\ell<m$.
Recall that for staircase polygons the figures of width one were an
exceptional case which did not obey the same reciprocity result as the
general case.  The staircase polygons with a hole of width one are also
an exceptional case, which we must exclude.  We thus impose the additional
condition $\ell-k>1$.  A figure with given $k$, $\ell$ and $m$ is specified by
the following dimensions, as shown in Figure~\ref{subfig:spshb},
\begin{enumerate}
\item heights $\Nb_j$ and $\Nt_j$ of the lower
and upper perimeter segments of the polygon, $1\le j\le m-1$,
\item interior heights $M_j$ to the left and right of, and within, the hole, 
$0\le j\le m$,
\item heights $\Hb_j$ and $\Ht_j$ of the lower and upper
perimeter segments of the hole, $k+1\le j\le \ell-1$,
\item interior heights $\Mb_j$ and $\Mt_j$ below and above
the hole, $k\le j\le \ell$.
\end{enumerate}

Three-choice polygons can be regarded as staircase polygons with a hole
which doesn't close.  The width $k$ has the same meaning as above, $\ell$
denotes the ultimate horizontal extent of the branch of the figure above
the hole, and $m$ denotes the ultimate horizontal extent of the branch
below the hole.
Note that $\ell\ge k$ and $m>k$.  
Again an exceptional case, $m=k+1$, must be excluded.  Hence we impose the
restriction $m>k+1$.
The labeling of the vertical dimensions follows, with
a few obvious modifications, the pattern of
staircase polygons with a staircase hole and is shown in
Figure~\ref{subfig:threechoice}.  In particular, the heights $M_j$ within
the hole are defined only for $j\le\min(\ell,m)$.  When $\ell=k$
the unknowns $\Mt_j$ and $\Ht_j$ do not appear.
This special case is treated separately.

As in the case of
column-convex polygons, the standard statistics are linear in these
heights.  The (half-)vertical perimeter for staircase polygons with a
staircase hole is given by
\begin{equation}
n = M_0+\sum_{j=1}^{m-1}\Nt_j +M_k+\sum_{j=k+1}^{\ell-1} \Ht_j
\label{vperimeterspsh}
\end{equation}
and the area is given by
\begin{equation}
a = M_0+\sum_{j=1}^{k-1}(M_j+\Nt_j) +\sum_{j=k+1}^\ell (\Nb_j+\Mb_j)
+\sum_{j=k}^{\ell-1} (\Mt_j+\Nt_j)+\sum_{j=\ell+1}^{m-1}(\Nb_j+M_j)+M_m.
\label{areaspsh}
\end{equation}
In what follows, we associate to the unknowns $\Nt_i$ (resp. $\Nb_i$,
$\Ht_i$, $\Hb_i$, $M_i$, $\Mt_i$, $\Mb_i$) the fugacities $\yt_i$
(resp. $ \yb_i, \wt_i, \wb_i, z_i, \zt_i, \zb_i$).
\begin{propo}\label{propo:spsh}
Let $E_{k,\ell,m}(\byt, \byb,\bwt, \bwb, \bz, \bzt, \bzb)$ be the
generating function for stair\-case
polygons with a staircase hole where $k$, $\ell$ and $m$ are the
widths defined above.
Then if $\ell-k>1$, the generating function
$E_{k,\ell,m}(\byt, \byb,\bwt, \bwb, \bz, \bzt, \bzb)$ is self-reciprocal,
$$E_{k,\ell,m}(1/\byt,1/ \byb,1/\bwt, 1/\bwb, 1/\bz, 1/\bzt, 1/\bzb)= 
\ \ \ \ \ \ \ \ \ \ \ \ \ \ \ \ \ \ \ \ \ \ \ \ \ \ \ \ \ \ \ \ \ \ \
\ \ \ \ \ \ \ \ \ \ \ \ \ $$
\begin{equation}\ \ \ \ \ \ \ \ \ \ \ \ -\frac{z_kz_\ell
\prod_{j=1}^{m-1}(\yt_j\yb_j)
\prod_{j=k+1}^{\ell-1}(\wt_j\wb_j)}
{\prod_{j=1}^{m-1}z_j
\prod_{j=k+1}^{\ell}\zb_j
\prod_{j=k}^{\ell-1}\zt_j}E_{k,\ell,m}(\byt, \byb,\bwt, \bwb, \bz,
\bzt, \bzb).
\label{selfspsh}
\end{equation}
\end{propo}
{\em Proof.} The linear homogeneous diophantine system is the union of
five sets of equations which we label $L_1$--$L_5$.  The regions to the left
and right of the hole give $L_1$ and $L_2$, the regions below and above
the hole give $L_3$ and $L_4$ and the inside of the hole gives $L_5$:
\begin{align} \label{Ldefsa}
L_1 &=\begin{cases}
  M_0-M_1-\Nb_1=0 \\
  M_j+\Nt_j-M_{j+1}-\Nb_{j+1}=0 \qquad \text{for $1\le j\le k-2$} \\
  M_{k-1}+\Nt_{k-1}-\Mt_k-M_k-\Mb_k-\Nb_k=0
\end{cases} \notag\displaybreak[0] \\[0.05in]
L_2 &=\begin{cases}
  \Mb_\ell+M_\ell+\Mt_\ell+\Nt_\ell-M_{\ell+1}-\Nb_{\ell+1}=0 \\
  M_j+\Nt_j-M_{j+1}-\Nb_{j+1}=0 \qquad \text{for $\ell+1\le j\le m-2$} \\
  M_{m-1}+\Nt_{m-1}-M_m=0
\end{cases} \notag\displaybreak[0] \\[0.05in]
L_3 &=\begin{cases}
  \Mb_k-\Mb_{k+1}-\Nb_{k+1}=0 \\
  \Mb_j+\Hb_j-\Mb_{j+1}-\Nb_{j+1}=0 \qquad \text{for $k+1\le j\le \ell-1$}
\end{cases} \notag\displaybreak[0] \\[0.05in]
L_4 &=\begin{cases}
  \Mt_j+\Nt_j-\Mt_{j+1}-\Ht_{j+1}=0 \qquad \text{for $k\le j\le \ell-2$} \\
  \Mt_{\ell-1}+\Nt_{\ell-1}-\Mt_\ell=0 
\end{cases} \notag\displaybreak[0] \\[0.05in]
L_5 &=\begin{cases}
  M_k-M_{k+1}-\Hb_{k+1}=0  \\
  M_j+\Ht_j-M_{j+1}-\Hb_{j+1}=0 \qquad \text{for $k+1\le j\le \ell-2$} \\
  M_{\ell-1}+\Ht_{\ell-1}-M_\ell=0.
\end{cases}
\end{align}
All heights of course are nonnegative.  Self-avoidance imposes the
additional constraint that the heights denoted $M_j$, $\Mt_j$ and
$\Mb_j$ be positive.  The set $\btau$, defined in step 2 above,
contains six unknowns whose constraints
are not independent which we eliminate as follows: $M_0$ using
the first equation of $L_1$,
$M_m$ using the last equation of $L_2$, $\Mb_k$ using the first equation of
$L_3$, $\Mt_\ell$ using the last equation of $L_4$, and $M_k$ and $M_\ell$ using
the first and last equations of $L_5$.  The resulting system is
\begin{align} \label{Lpdefs}
L'_1 &=\begin{cases}
  M_j+\Nt_j-M_{j+1}-\Nb_{j+1}=0 \qquad \text{for $1\le j\le k-2$} \\
  M_{k-1}+\Nt_{k-1}-\Mt_k-M_{k+1}-\Hb_{k+1}-\Mb_{k+1}-\Nb_{k+1}-\Nb_k=0 &
\end{cases} \notag\displaybreak[0] \\[0.05in]
L'_2 &=\begin{cases}
  \Mb_\ell+M_{\ell-1}+\Ht_{\ell-1}+\Mt_{\ell-1}+\Nt_{\ell-1}+\Nt_\ell-M_{\ell+1}-\Nb_{\ell+1}=0 \\
  M_j+\Nt_j-M_{j+1}-\Nb_{j+1}=0 \qquad \text{for $\ell+1\le j\le m-2$}
\end{cases} \notag\displaybreak[0] \\[0.05in]
L'_3 &=\begin{cases}
  \Mb_j+\Hb_j-\Mb_{j+1}-\Nb_{j+1}=0 \qquad \text{for $k+1\le j\le \ell-1$}
\end{cases} \notag\displaybreak[0] \\[0.05in]
L'_4 &=\begin{cases}
  \Mt_j+\Nt_j-\Mt_{j+1}-\Ht_{j+1}=0 \qquad \text{for $k\le j\le \ell-2$}
\end{cases} \notag\displaybreak[0] \\[0.05in]
L'_5 &=\begin{cases}
  M_j+\Ht_j-M_{j+1}-\Hb_{j+1}=0 \qquad \text{for $k+1\le j\le \ell-2$.}
\end{cases}
\end{align}

The substitutions
$\Nt_j$, $\Nb_j$, $\Ht_j$, $\Hb_j=1$ and $M_j$, $\Mt_j$, $\Mb_j=-1$ solve this
new system of equations.  One should note that when $k=1$ or $m-\ell=1$ the
system is somewhat modified, but one may check that the solution still holds.
Therefore we may apply Corollary~\ref{coro:selfb} to obtain a self-reciprocity
condition on the generating function for the solutions of $\bigcup_j L'_j$
subject to the positivity constraints on the heights.  
%
%
Making appropriate substitutions
to restore the unknowns in set $\btau$, and using Lemma~\ref{lemma:stability}
we obtain eqn.~(\ref{selfspsh}).
\qed

\medskip
We now treat three-choice polygons.
\begin{propo}\label{propo:three-choice}
 Let $E_{k,\ell,m}(\byt,\byb,\bwt, \bwb, \bz,\bzt,\bzb)$ be the generating function for three-choice
polygons where $k$, $\ell$ and $m$ are the widths defined above.
Then if $m-k>1$, the generating function
$E_{k,\ell,m}(\byt,\byb,\bwt, \bwb, \bz,\bzt,\bzb)$ satisfies a self-reciprocity condition which, when
$\ell=k$, takes the form
$$E_{k,k,m}(1/\byt,1/\byb,1/\bwt, 1/\bwb, 1/\bz,1/\bzt,1/\bzb)=\ \ \ \ \
\ \ \ \ \ \ \ \ \ \ \ \ \ \ \ \ \ \ \ \ \ \ \ \ \ \ \ \ \ \ $$
\begin{equation}\ \ \ \ \ \ \ \ \ \ \ \ -\frac{\prod_{j=1}^{k-1}\yt_j \prod_{j=1}^{m-1}\yb_j
\prod_{j=k+1}^{m-1}\wb_j}
{\prod_{j=1}^{k}z_j\prod_{j=k+1}^{m-1}\zb_j}E_{k,k,m}(\byt,\byb,\bwt, \bwb, \bz,\bzt,\bzb),
\label{selfthreechoiceff}
\end{equation}
and, when $\ell>k$, takes the form
$$E_{k,\ell,m}(1/\byt,1/\byb,1/\bwt, 1/\bwb, 1/\bz,1/\bzt,1/\bzb)=\ \ \ \
\ \ \ \ \ \ \ \ \ \ \ \ \ \ \ \ \ \ \ \ \ \ \ \ \ \ \ \ \ \ \ $$
\begin{equation}
\ \ \ \ \ \ \ \ \ \ \ -\frac{z_k\prod_{j=1}^{\ell-1}\yt_j\prod_{j=1}^{m-1}\yb_j
\prod_{j=k+1}^{\ell-1}\wt_j\prod_{j=k+1}^{m-1}\wb_j}
{\prod_{j=1}^{\min(\ell,m-1)}z_j\prod_{j=k}^{\ell-1}\zt_j\prod_{j=k+1}^{m-1}\zb_j}E_{k,\ell,m}(\byt,\byb,\bwt, \bwb, \bz,\bzt,\bzb).
\label{selfthreechoicefg}
\end{equation}
\end{propo}
{\em Proof.} It is simpler to treat the two cases $\ell=k$ and $\ell>k$ separately.
The proofs follow very closely that of Proposition~\ref{propo:spsh}. \qed

As for column-convex polygons, the two propositions above may be
extended to other statistics.

\begin{coro} \label{coro:linearb}
Let ${\cal P}$ be either of the sets staircase polygons with a staircase
hole or three-choice polygons.  Let ${\cal P}_{k,\ell,m}$ be the subset
of figures in ${\cal P}$ with the widths $k$, $\ell$ and $m$ defined
as above.  Then the generating function in ${\cal P}_{k,\ell,m}$ according
to any statistics linear in the quantities  $\bNt, \bNb,$ 
$\bHt$, $\bHb$, $\bM$, $\bMt$, $\bMb$,
 is a
self-reciprocal rational function 
(assuming $\ell-k>1$ for staircase
polygons with a staircase 
hole and $m-k >1$ for three-choice polygons).
\end{coro}

The half-horizontal perimeter for either of the sets $\cal P$ is
given by $m+\ell-k$.  Using this in combination with
Corollary~\ref{coro:linearb},~\Ref{vperimeterspsh} and~\Ref{areaspsh}, we
obtain the inversion relations specialized to horizontal and vertical
perimeter, and area, which are listed in Table~\ref{table:summary}.
The exceptional cases ($\ell-k=1$ and  $m-k =1$ respectively) can be computed
explicitly by the methods of~\cite{bousc}.
\section{Discussion}
Each of the methods we have discussed for obtaining reciprocity or
inversion relations has its own particular uses.  For example, the
method of Stroganov is suitable for lattice models in
statistical mechanics which are characterized by a family of commuting
transfer matrices.  The Temperley methodology is mainly applicable
to families of polygons that are column-convex or nearly so.  Stanley's
method for obtaining reciprocity results apply to any problem defined by a
system of linear homogeneous diophantine (LHD) equations, but the solutions to
this system
must be constrained by a system of simple inequalities of a certain form.

It is probable that for many lattice models in statistical mechanics the
low temperature expansion can be framed as an LHD-system.  However,
most are likely
to require more general types of constraints than the simple inequalities of
the directed polyomino problems we have considered.  Likewise, the
non-directed polygon problems that we have successfully treated using
the Temperley methodology can be recast as LHD-systems with more complex
constraints.  How to handle such constraints is a worthy problem for
future investigation.

In recent work~\cite{beo} this statistical mechanical language
has been adapted for the enumeration of lattice paths,
and may apply to polyomino problems as well.  It is intriguing
to speculate that the inversion relations found here may be connected
with this approach.

We have not searched for inversion relations for any polyomino problem
in variables other than the natural variables for the problem.  Yet the
example of the Potts model demonstrates that such inversion relations
may exist.  It is also possible that symmetries in addition to the ones
presented here can be found for some problems.  It is our hope that such
additional symmetries might lead to the solution of currently intractable
problems.

For the moment, we remark that the search for inversion and symmetry
relations appears to provide a new method to tackle certain combinatorial
problems.  The degree of applicability of this method is still unclear.

\section*{Acknowledgments}
We have benefited from conversations with George Andrews, Richard Brak,
Jean-Marie Maillard, Paul Pearce and Markus V{\"o}ge and from
correspondence with Jean-Marc F{\'e}dou and Richard Stanley.  We
thank Iwan Jensen for providing us with his series data for
staircase polygons with a staircase hole.
AJG and WPO acknowledge support from the Australian Research Council.

\end{document}